\RequirePackage{fix-cm}
\documentclass[natbib,runningheads]{svjour3}
\smartqed  % flush right qed marks, e.g. at end of proof
\usepackage{graphicx}
\usepackage{bm}
\usepackage{natbib}
\usepackage{amssymb}
\usepackage{theorem}
\usepackage{xspace}
\usepackage{amsfonts}
\theoremstyle{plain}
\usepackage{graphicx}
\usepackage{amsmath}
\usepackage{array}
\usepackage{float}
\usepackage{color}
 \usepackage{subfigure}
\usepackage{epsfig}
\usepackage{hhline}
\usepackage[utf8]{inputenc}    % pour les accents

 \numberwithin{equation}{section} 
\numberwithin{figure}{section}
\def\EL{\textsc{E\!L}}
\def\Mx{\textsc{M\!x}}

\def\eeta{\textrm{\mathversion{bold}$\mathbf{\eta}$\mathversion{normal}}}

\def\eb{\textrm{\mathversion{bold}$\mathbf{\beta}$\mathversion{normal}}}  \def\del{\textrm{\mathversion{bold}$\mathbf{\gamma}$\mathversion{normal}}}

\def\ebo{\textrm{\mathversion{bold}$\mathbf{\beta^0}$\mathversion{normal}}}
\def\el{\textrm{\mathversion{bold}$\mathbf{\lambda}$\mathversion{normal}}}

\def\epsi{\textrm{\mathversion{bold}$\mathbf{\psi}$\mathversion{normal}}}

\def\R{I\!\!R}
\def\eR{I\!\!R}
\def\eE{I\!\!E}
\def\eP{I\!\!P}
\def\e1{1\!\!1}
\def\e1{1\!\!1}

\def\eg{\mathbf{{G}}}

\def\eV{\mathbf{{V}}}
\def\ebb{\mathbf{{B}}}
\def\ek{\mathbf{{K}}}
\def\er{\mathbf{{R}}}
\def\es{\mathbf{{S}}}

\def\hh{\hspace*{0.3cm}}

\def\ev{\mathbf{{V}}}

\def\em{\mathbf{{M}}}

\def\ew{\mathbf{{w}}}
\def\en{\mathbf{{N}}}
\def\ez{\mathbf{{z}}}

\def\tr{\text{tr}}
\newcommand{\N}{\mathbb{N}}
\def\XX{\textrm{\mathversion{bold}$\mathbf{X}$\mathversion{normal}}}

\def\ell{\textrm{\mathversion{bold}$\mathbf{\overset{.}{\lambda}}$\mathversion{normal}}}
\newcommand{\Var}{\mathbb{V}\mbox{ar}\,}

\newtheorem{Remark}{Remark} %[section]
\newtheorem{Corollary}{Corollary}
\begin{document}
\title{Empirical likelihood test for high-dimensional two-sample model}

\author{Gabriela Ciuperca         \and
        Zahraa Salloum 
        }

%\authorrunning{Short form of author list} % if too long for running head

\institute{ G.Ciuperca - Z. Salloum \\
Universit\'e de Lyon, Universit\'e Lyon 1, CNRS, UMR 5208, Institut Camille Jordan, \at
             Bat.  Braconnier, 43, blvd du 11 novembre 1918, F - 69622 Villeurbanne Cedex, France \\
                         \email{Gabriela.Ciuperca@univ-lyon1.fr}           \\ \\
            % \emph{Present address:} of F. Author  %  if needed
           %\and
       %   Universit\'e de Lyon, Universit\'e Lyon 1, CNRS, UMR 5208, Institut Camille Jordan, \at
         %    Bat.  Braconnier, 43, blvd du 11 novembre 1918, F - 69622 Villeurbanne Cedex, France \\
       %    Tel.: 33(0)4.72.44.62.19\\
        %   Fax: 33(0)4.72.43.16.87\\
    Z. Salloum \\
            \email{salloum@math.univ-lyon1.fr}      
}
\maketitle
\begin{abstract}
A non parametric method based on the empirical likelihood is proposed for detecting the change   in the coefficients of   high-dimensional linear model where the number of model variables may increase as the sample size increases.  This amounts to testing the null hypothesis of no change against the alternative of one change in the regression coefficients. Based on the theoretical asymptotic behaviour of the empirical likelihood ratio statistic, we propose, for a fixed design, a simpler test statistic, easier to use in practice.    The asymptotic normality of the proposed test statistic under the null hypothesis is proved, a result which is different from the $\chi^2$ law for a model with a fixed variable number. Under alternative hypothesis, the test statistic diverges. We can then find  the asymptotic confidence region for the difference of parameters of the two phases. Some Monte-Carlo simulations study the behaviour of the proposed test statistic.
\end{abstract}
%%%%%%%%%%%%%%%

%\textit{AMS 2000 subject classifications} : Primary 62F03, 62G10; secondary 62G10, 62F05.

%%%%%%%%%%%%%%%%%%%%%%

\section{Introduction} 
\label{sec1}
\hh The technology development and  fast numerical techniques make possible to consider and study statistical models with a large number of variables.   High-dimensional model refers to a model  whose the number  $p$ of explanatory variables increases to infinity as the number $n$ of observations converges to infinity. When $p$ diverges, traditional statistical methods may not work with this kind of growth dimensionality.  \\
\hh Most of the literature works on high-dimensional model utilize the LASSO (Least Absolute Shrinkage and Selection Operator) type methods, in order to automatically select the significant  variables. The principle of these methods, introduced by \cite{Tibshirani:96},  is to optimize a penalized process, more precisely, a process with a $L_1$-type penalty. If the model contains outliers, the parameter estimators by the least squares method with LASSO penalty have a large error. An alternative method is then the penalized quantile method. Thereby, \cite{Dicker:Huang:Lin:13} consider a quantile model with seamless-$L_0$ penalty when the number $p$ of explanatory variables  is such that $p \rightarrow \infty$, $p /n \rightarrow 0$ as $n \rightarrow \infty$. For a general quantile regression, \cite{Wu:Liu:09} propose the SCAD penalty, while, in  \cite{Zou:Yuan:08}, a composite quantile regression is considered with an adaptive LASSO penalty. The case $p \rightarrow\infty$ is also considered in \cite{Fan:Peng:04} for a non-concave penalized likelihood method, when $p^5/n \rightarrow \infty$. Concerning the group selection methods for high-dimensional models, the readers find in \cite{Huang:Breheny:Ma:12} a review of  methods. \\
All these methods are based first on the principle of selecting (automatically) the significant variables. Then,  the dependent variable is  modeled only as a  function of the significant variables, in order to have more accurate parameter estimators and a  better adjustment for the dependent variable. \\
\hh If the goal is to have the most accurate prediction and also robust, in the case of a model with outliers, one possibility is to consider the empirical likelihood (EL)  method. But, for this type of method, in literature,  most papers are devoted  to the case of  fixed $p$. For a high-dimensional linear regression model, we can refer first to the paper of \cite{Guo:Zou:Wang:Chen:13}, when the design is deterministic. High-dimensional data are also studied by  \cite{Liu:Zou:Wang:13}, where EL method is considered for a sequence of i.i.d. random vectors with dimension $p$,  when $p \rightarrow \infty$ as $n \rightarrow \infty$. \\
\hh In this paper, we are interested by a change-point model, that is, a model which changes at some moment. The number $p$ of explanatory variables varies with the number $n$ of observations and $p$ can converge to infinity if $n \rightarrow \infty$. \\
\hh Since statistical techniques in high-dimension are fairly recent, there are not many papers in literature that address the change-point problem in a high-dimensional model. \cite{Lung-Yut-Fong:Levy:Cappe:12} propose an approach for detection of a change-point in high-volume network traffic. The asymptotic distribution of the  test statistic proposed in \cite{Lung-Yut-Fong:Levy:Cappe:12}, under the null hypothesis that there is no change-point, is the \textit{argsup} of a Brownian Bridge. There are some papers where  LASSO type methods are used.   \cite{Lee:Seo:Shin:12} consider a possible change-point in a high-dimensional regression with Gaussian errors.  The main result of the article is to show that the sparsity property is maintained, even if there is a change in the model. There is no hypothesis test to decide the presence or absence of change in model. In \cite{Ciuperca:13a}, LASSO-type and adaptive LASSO estimators  are studied, while   in \cite{Ciuperca:13b} quantile model with SCAD penalty is considered. These last two papers consider models with $p$ fixed.  In order to choose the change-point number, a model selection criterion is also proposed by \cite{Ciuperca:13a}.  \\
\hh To the authors' knowledge, the  EL technique has not yet been addressed in a high-dimensional two-sample model,  that makes the interest of this work.  We study the asymptotic behaviour of the empirical likelihood ratio test statistic  when the  design is deterministic. \ \\  
%%%%%%%%%%%%%%%%%%%%%%

We consider a first linear model: 
\begin{equation}
\label{eq1}
Y_i=\XX_i^t\eb+\varepsilon_i, \qquad i=1, \cdots, n.
\end{equation}
Consider now  a second linear model which changes at observation $k$. It is called two-phase model, or model with one change-point: 
\begin{equation}
\label{eq2}
Y_i=\left\{
\begin{array}{ccl}
\XX_i^t\eb+\varepsilon_i, \qquad 1 \leq i \leq k, \\
\XX_i^t\eb_2+\varepsilon_i, \qquad k < i \leq n,
\end{array}
\right.
\end{equation}
where $\XX_i$ is a $p \times 1$ vector of $p$  explanatory variables, $\eb$ and $\eb_2$ are $p\times 1$ vectors of unknown parameters and $\varepsilon_i$ designates the model error.
The parameter $\eb$ of the first phase of (\ref{eq2}) coincides with that of (\ref{eq1}). For  models (\ref{eq1}) and (\ref{eq2}), $Y_i$ is observation $i$ of the response variable. 
The errors  $\varepsilon_i$ are supposed  independent identically distributed (i.i.d), with mean zero and finite variance $\sigma^2$.\\

 We assume that the number $p$ of explanatory variables  $\XX_i$ depends on the sample size $n$: $p=p_n$, such that $p_n \rightarrow \infty$ as $n \rightarrow \infty$. The change-point $k$ of (\ref{eq2}) also depends on $n$. The change in model (\ref{eq2}) takes place far enough from the first observation and sufficiently previous to the last observation. So, we suppose that $\lim _{n \rightarrow \infty} k/n \in (0,1)$.
\\ 

In this paper, for given  $k$, we use the empirical likelihood method to construct the confidence region for $\eb-\eb_2$, or equivalent to test the null hypothesis of no change in  model (\ref{eq2}). Under null hypothesis, the model has the form (\ref{eq1}), that is 
\begin{equation}
\label{H0}
H_0 : \eb_2=\eb.
\end{equation}
The alternative hypothesis assumes that one change occurs in the regression parameters, that is 
\begin{equation}
\label{H1}
H_1 : \eb_2\neq \eb.
\end{equation}

The paper is organized as follows. In Section \ref{sec2} we first present the EL method for the   two-sample model. Some notations  used throughout the paper are defined and  needed assumptions for the theoretical study are also announced. In Section \ref{sec3}, we construct an empirical likelihood ratio test statistic and we study its asymptotic behaviour. The asymptotic distribution under $H_0$ of the test statistic is obtained, while, under $H_1$, this statistic diverges. Next, in Section \ref{sec4}, we analyse  the coverage accuracy and the empirical power by means of simulations, which confirm the performance of proposed test.  A new critical value is also proposed in order to improve the coverage rate. The proofs of the main results are given in Appendix (Section \ref{sec5}) followed by some Lemmas and their proofs. 
%%%%%%%%%%%%%%%%%%%%%%%%%%%%%%%%%%%%%%%%%%%%%%%%%%%%%%%%%%%
%%%%%%%%%%%%%%%%%%%%%%%%%%%%%%%%%%%%%%%%%%%%%%%%%%%%%%%%%%%
\section{Preliminares}
\label{sec2}
\hh In this section, we introduce the EL method for the  two-sample model. Notations and assumptions are also given. \\

Under null hypothesis $H_0$, that is  model (\ref{eq1}), let $\ebo$ denote the true value   of the parameter $\eb$. Under alternative hypothesis $H_1$, that is  model (\ref{eq2}), the  true values  of $\eb$, $\eb_2$, respectively, are  $\eb^0$, $\eb^0_2$.\\
In order to define the profile empirical likelihood (under $H_0$ and under $H_1$), we introduce the following  random $p$-vector, for all $\eb \in \R^p$ and $i=1, \cdots, n$:
\begin{equation*}
\label{eq3}
\ez_i(\eb) \equiv \XX_i(Y_i-\XX_i^t\eb).
\end{equation*}
Consider also the vector \[\ez_i^0 \equiv \XX_i \varepsilon_i.\]
We remark that, under the hypothesis $H_0$, we have $\ez_i^0 = \ez_i(\ebo)$, for all $i=1, \cdots, n$ and $\eE[\ez_i^0]=\textbf{0}_p$. On the other hand, for fixed design $(\XX_i)_{1 \leq i \leq n}$, the random variables  $\ez_i^0$ are independent but not identically distributed. \\
On the change-point, we consider the notation $\theta_{nk}=k/n$. Thus, in view of the remark made in Introduction, we assume that $\theta_{nk}\rightarrow \theta^0 \in (0,1)$ as $n \rightarrow \infty$.\\ 

For  the dependent variable $Y_i$ of  model (\ref{eq2}), let us consider the probability to observe the value $y_i$ (respectively $y_j$) : $q_i\equiv \eP [Y_i=y_i]$, for $i=1, \ldots, k$ and $q_j \equiv \eP[Y_j=y_j]$, for $j=k+1, \cdots,n$. Obviously, these probabilities satisfy the relations $\sum  _ {i=1} ^ {k} q_i = 1$ and $\sum   _{j =k+1} ^n q_j = 1$. Corresponding to these probabilities, we define the probability vectors $(q_1, \cdots, q_{k})$ and $(q_{k+1}, \cdots, q_{n})$.
\\  
%%%%%%%%%%%%%%%%%%%%%%%%%%%%%%%%%%%%%%%%%%%%%%%%%%%%%
 Under  hypothesis $H_0$ given by (\ref{H0}), the profile empirical likelihood for $\eb$ is   
\[
{\cal R}_{nk}(\eb) \equiv \sup   _{(q_1,\cdots,q_k)}
  \sup   _{(q_{k+1},\cdots,q_n)}  \big\{  \prod    _{i=1} ^ {k}q_i   \prod   _{j =k+1} ^n q_j ;   \sum    _{i=1} ^ {k} q_i = 1,  \sum    _{j =k+1} ^n q_j =1,
  \sum    _{i=1} ^ {k}q_i \ez_i( \eb )= \sum _{j =k+1} ^n
q_j \ez_j( \eb)=\textbf{0}_p \big \},
\]
with $\textbf{0}_p$ the $p$-vector with all components zero. \\
Similarly, under hypothesis $H_1$ given by (\ref{H1}), the profile empirical likelihood is 
\[
{\cal R}_{nk}(\eb, \eb_2)  \equiv \sup  _{(q_1,\cdots,q_{k})} \sup  _{(q_{k+1},\cdots,q_n)} 
  \big\{  \prod  _ {i=1} ^{ k} q_i   \prod  _ {j =k+1} ^{n } q_j ;   \sum  _ {i=1} ^ {k} q_i =1,     \sum _ {j =k+1} ^{n } q_j =1,  \sum  _ {i=1} ^{ k}q_i \ez_i( \eb )= \textbf{0}_p, \sum  _{j =k+1} ^ n q_j \ez_j( \eb_2)=\textbf{0}_p \big\}.
\]
Then, using an idea similar to the maximum likelihood test for testing $H_0$ against $H_1$, we consider the profile empirical likelihood ratio
$
{{\cal R}_{nk}(\eb)}/{{\cal R}_{nk}(\eb, \eb_2)}$.\\
 Since
$
{\cal R}_{nk}(\eb, \eb_2) = k^{-k} (n-k)^{-(n-k)}$,
we have that the corresponding empirical log-likelihood ratio  is 
\[
 -2 \sup   _{(q_1,\cdots,q_{k})} 
  \sup   _{(q_{k+1},\cdots,q_n)} \big \{ \sum   _ {i=1} ^ {k} \log (kq_i ) + \sum  _{j =k+1} ^ n \log ((n-{k}) q_j);   \sum  _ {i=1} ^ {k} q_i =   \sum _{j =k+1} ^n q_j =1, \sum    _{i=1} ^ {k}q_i \ez_i( \eb )= \sum _{j =k+1} ^n
q_j \ez_j( \eb)=\textbf{0}_p
 \big \}.
\]
Applying the Lagrange multiplier method, the optimal probabilities $q_i$ and $q_j$ are
\begin{equation}
\label{qiqj}
q_i= \frac{1}{k +n \el_1^t \ez_i( \eb ) }, \qquad q_j= \frac{1}{n-k- n  \el_2^t \ez_j( \eb ) },
\end{equation}
where $\el_1, \el_2 \in \eR^p$ are the Lagrange multipliers.
Consequently, the corresponding empirical log-likelihood function can be written as
\begin{equation}
\label{eq5} 
 2 \sum   _ {i=1} ^ {k}\log
\big({1+\frac{n}{k}\el_1^t \ez_i(\eb)}\big)+ 2\sum  _{j =k+1} ^ n \log \big({1-\frac{n}{n-k}\el_2^t \ez_j(\eb)}\big).
\end{equation}
Taking into account relation (\ref{qiqj}), the derivative with respect to $\eb$ of (\ref{eq5}) is $\sum^k_{i=1} q_i \XX_i \XX_i^t \el_1- \sum^n_{j=k+1} q_j \XX_j \XX_j^t \el_2=\textbf{0}_p$. We can apply  Lemma 4 of \cite{Guo:Zou:Wang:Chen:13}  on each phase of model, that implies that  $\| \el_1\|=O_{\eP}(p^{1/2} k^{-1/2})$ and $\| \el_2\|=O_{\eP}(p^{1/2} (n-k)^{-1/2})$. Then, the probabilities $q_i$ and $q_j$ of (\ref{qiqj}) are approximatively $k^{-1}$ and $(n-k)^{-1}$, respectively. Thus, we can restrict $\el_1$ and $\el_2$ such that 
\begin{equation}
\label{VkVn}
k^{-1}\sum^k_{i=1}  \XX_i \XX_i^t \el_1=(n-k)^{-1} \sum^n_{j=k+1} \XX_j \XX_j^t \el_2. 
\end{equation}
If the symmetric matrices $k^{-1}\sum^k_{i=1}  \XX_i \XX_i^t$ and $(n-k)^{-1} \sum^n_{j=k+1} \XX_j \XX_j^t$ converge, as $n \rightarrow \infty$, to two strictly positive definite matrices, then the relation (\ref{VkVn}) can be written $\el_1=\big( k^{-1}\sum^k_{i=1}  \XX_i \XX_i^t \big)^{-1} \big( (n-k)^{-1} \sum^n_{j=k+1} \XX_j \XX_j^t \big)  \el_2$. Noting by $\widetilde \el_2 \equiv \big( k^{-1}\sum^k_{i=1}  \XX_i \XX_i^t \big)^{-1} \big( (n-k)^{-1} \sum^n_{j=k+1} \XX_j \XX_j^t \big)  \el_2 $, we have the new Lagrange multipliers such that $\el_1= \widetilde \el_2$. For the sake of readability, we denote $\widetilde \el_2$ by $\el_2$.\\

With this remark, we will restrict the study to a particular case, when $\el_1=\el_2=\el$. Considering this constraint, instead of statistic (\ref{eq5}) we consider the following particular empirical likelihood ratio (ELR) statistic
\begin{equation}
\label{eq6}
\EL_{nk}(\eb) \equiv 2 \sum  _ {i=1} ^ {k}\log
\big(1+\frac{n}{k} \el^t\ez_i(\eb)\big)+ 2\sum  _{j =k+1} ^ n \log
\big(1-\frac{n}{n-k}\el^t\ez_j(\eb)\big),
\end{equation}
where the Lagrange multiplier $\el \in \eR^p$ satisfies
\begin{equation} 
\label{eq120}
 \sum  _ {i=1} ^ {k}  \frac{\ez_i(\eb)}{k/n +\el^t \ez_i(\eb)} - \sum  _{j =k+1} ^ n \frac{\ez_j(\eb)}{1- k/n  -\el^t \ez_j(\eb)} = \textbf{0}_p.
\end{equation}
%\label{Sub Notations}
\subsection{Notations}
%\label{Sec3}
We provide a brief summary of notations used in the paper. \\

For exposition convenience, we define some general notation. All vectors are column and $\textbf{v}^t$ denotes the transposed of \textbf{v}. All vectors and matrices are in bold.  For a vector \textbf{v}, by $\|\textbf{v}\|$ we denote its Euclidean norm and by $\|\textbf{v}\|_1$ its $L_1$-norm. For a symmetric p-square matrix $\textbf{A}=(a_{ij})$, let us denote by $\gamma_1(\textbf{A}) \geq \gamma_2(\textbf{A}) \geq \ldots \geq \gamma_p(\textbf{A})$ the eigenvalues  and $\tr(\textbf{A})$ as the trace operator of the matrix $\textbf{A}$. Consider also the following notation  $\Mx(\textbf{A})= \max _{1 \leq i,j \leq p} | a_{ij}|$. We  denote by $\| \textbf{A} \|_1= \max_{j=1,\cdots,p} (\sum _{i=1}^p |a_{ij}|)$, the subordinate norm to the vector norm $\| . \|_1$.\\
 All throughout the paper, C denotes a generic constant which may be different from line to line and even from formula to formula and whose value is not of interest.\\
  Moreover, $\textbf{0}_p$ denote the $p-$vector with all components zero.  \\

At the beginning of this section, the notation   $\theta_{nk} \equiv k/n$ was introduced. To simplify notations, we will use the notation $\theta$ instead $\theta_{nk}$. \\ 

For $\ebo$, the true value of the parameter $\eb$ on the phase $1, \cdots, k$, and the test value under $H_0$, we define the following  $p$-square matrix
\begin{equation}
\label{Sb0}
\es_n(\ebo) \equiv \frac{1}{n \theta ^2 } \sum  _ {i=1} ^ {k}  \ez_i(\eb^0) \ez_i^{t}(\eb^0) + \frac{1}{n (1-\theta)^2} \sum _{j =k+1} ^ n \ez_j(\eb^0) \ez_j^{t}(\eb^0)
\end{equation}
and the following  $p$-vector 
\begin{equation}
\label{psib}
\epsi_n(\ebo) \equiv  \frac{1}{n \theta} \sum  _ {i=1} ^ {k}  \ez_i(\eb^0) - \frac{1}{n (1-\theta )} \sum _{j =k+1} ^ n \ez_j(\eb^0).
\end{equation}
Under null hypothesis, for the true value $\ebo$ of $\eb$, the mean of the random matrix $\es_n(\ebo)$ is the following $p$-square matrix 
\begin{equation}
\label{eq4}
\eV_n ^0 \equiv \frac{1}{n \theta ^2 } \sum  _ {i=1} ^ {k} \ev_{(i)}^0 + \frac{1}{n (1-\theta )^2} \sum _{j =k+1} ^ n \ev_{(j)}^0 ,
\end{equation}
where, for $i=1,\ldots, n$ 
\begin{equation}
\label{var}
\ev_{(i)}^0\equiv \Var(\ez_i^0)= \sigma^2 \XX_i \XX_i^t. 
\end{equation}
For $ i = 1, \cdots, n$, let us also consider the following random vector \[\ew_i^0 \equiv (\eV_n^0)^{-1/2} \ez_i(\ebo).\]
Corresponding to the components of $\ew_i^0=(w_{i,1}^0 , \ldots, w_{i,p}^0 )$,  we consider for $ i =1, \ldots,n$, for $r \in \N^*$, $t_1, \cdots, t_r \in \{1, \cdots , p\}$,   the following scalar
%We need also the some additional notation
\begin{equation}
\label{eq8}
\alpha^{t_1 t_2 \cdots t_r}\equiv \frac{1}{n\theta^r}  \sum  _ {i=1} ^ {k} \eE\big[ w_{i,t_1}^0 w_{i,t_2}^0 \cdots w_{i,t_r}^0 \big] +\frac{1}{ n(\theta-1)^r} \sum _{j =k+1} ^ n \eE\big[ w_{j,t_1}^0 w_{j,t_2}^0 \cdots w_{j,t_r}^0 \big]
\end{equation}
and the following  random variable
\begin{equation}
\label{eq9}
\omega^{t_1 t_2 \cdots t_r} \equiv \frac{1}{n \theta^r} \sum  _ {i=1} ^ {k} w_{i,t_1}^0 w_{i,t_2}^0 \cdots w_{i,t_r}^0 +\frac{1}{n (\theta-1)^r} \sum _{j =k+1} ^ n w_{j,t_1}^0 w_{j,t_2}^0 \cdots w_{j,t_r}^0 -\alpha^{t_1 t_2 \cdots t_r},
\end{equation}
where  $w_{i,t_r}^0 $ is the $r$-th component of $\ew_i^0$. In particular, for all $t_1, t_2 \in \{ 1, \cdots , p\}$,  $\alpha^{t_1}=0$, $\alpha^{t_1 t_2}= \delta ^{t_1t_2}$ is the Kronecker delta, that is $\alpha^{t_1t_2}=1$ if $t_1=t_2$, and 0 otherwise.

%%%%%%%%%%%%%%%%%%%%%%%%%%%%%%%%%%%%%%%%%%%%%%%%%%%%%%%%%%%%%%%%%%%%%%%%%%
%%%%%%%%%%%%%%%%%%%%%%%%%%%%%%%%%%%%%%%%%%%%%%%%%%%%%%%%%%%%%%%%%%%%%%%%%%%%%%%%ùù
%%%%%%%%%%%%%%%%%%%%%%%%%%%%%%%%%%%%%%%%%%%%%%%%%%%%%%%%%%%%%%%%%%%%%%%%%%%%%%%%%%%%%%%%%%%
\subsection{Assumptions}
\hh We now state the assumptions on the design, on the errors, on the number $p$ of the explanatory variables and on the change-point location. These assumptions are needed in order to keep the properties obtained for EL statistic in a high-dimensional model, without  change-point. \\ 
For assumptions (A3)-(A6) the constant $q$ is such that $q \geq 4$.\\ 

\noindent {\bf (A1)} There exist positive constants $C_0, C_1 >0$, such that $0 < C_0 < \inf _n  \gamma_{ { 1}} (\eV_n^0) \leq \sup_n  \gamma_1 (\eV_n^0) < C_1 < \infty$.
 \\ 
{\bf (A2)} $\eE(\varepsilon_1^4) < C_2$  for some $C_2>0$ and for all $n$.  \\  
{\bf (A3)} $p^{-1} \sum ^{p}_{s=1} |X_{i,s} |^q <C_3, 1\leq i \leq n$, for some $C_3 >0$, and $q \geq 4$;\\
{\bf (A4)} $\eE|\epsilon_1|^{2q} < C_4$, for some $C_4 >0$. \\
{\bf (A5)} $ p \ k^{(2-q)/(2q)} \rightarrow 0$ and $ p \ (n-k)^{(2-q)/(2q)} \rightarrow 0$, as $n \rightarrow \infty$. \\
{\bf (A6)} $ p^{2+4/q}\ k^{-1} \rightarrow 0$ and $ p^{2+4/q}\ (n-k)^{-1} \rightarrow 0$, as $n \rightarrow \infty$.\\
{\bf (A7)} $ \sum ^{p}_{r,s=1} \alpha^{rrss} =O(p^2)$.\\
{\bf (A8)} $ \sum ^{p}_{r,s,u=1} \alpha^{rsu} \alpha^{rsu} =O(p^{5/2})$ and $ \sum ^{p}_{r,s,u=1} \alpha^{rss} \alpha^{suu} =O(p^{5/2})$.\\ 
{\bf (A9)} For all $ i=1,\cdots, n$, for $l \in \N^*$ ,  $j_1,\cdots, j_l \in \{1, \cdots, p\}$ , and whenever $\sum_{i=1} ^l d_i \leq 6$, there exists a positive absolute constant $C_5< \infty$, then
$
\eE( w_{i,j_1} ^{d_1} \cdots  w_{i,j_l} ^{d_l}) \leq C_5$.\\

Assumptions (A3) and (A6) guarantee that the eigenvalues of $\es_n^0$ are close to those of $\eV_n^0$  (see Lemma \ref{corolary2}).  Assumption (A1) implies that  $\eV_n^0$ is uniformly nonsingular and bounded, for large $n$. Then, for $n$ large enough,  with probability close to one, $\es_n^0$ is nonsingular and $0< C_0 < \gamma_p(\es_n^0) \leq \gamma_1(\es_n^0) <C_1 < \infty$. Assumption (A3) is also assumed by \cite{Guo:Zou:Wang:Chen:13}, \cite{Hjort:Mckeague:Vankeilegom:09}, \cite{Liu:Zou:Wang:13} for high-dimensional model without change-point. Assumption (A4) together with (A3) and (A6) imply $\sup_{1\leq i \leq n} |\el^t \ez_i^0 |=o_p(1)$, which leads to Taylor expansions of (\ref{eq6}) and (\ref{eq120}) (see Lemma \ref{lemma5}). Assumptions (A1), (A2),  (A4) are also used by \cite{Guo:Zou:Wang:Chen:13} for linear models without change-point with random design.  Same assumption (A1) is requested in  \cite{Zi:Zou:Liu:12} for a two-sample model with fixed $p$. Assumptions (A5)-(A9) are also assumed by \cite{Guo:Zou:Wang:Chen:13}, \cite{Liu:Zou:Wang:13}, in order to have  for the asymptotic normality of the ELR statistic.

%%%%%%%%%%%%%%%%%%%%%%%%%%%%%%%%%%%%%%%%%%%ù
%%%%%%%%%%%%%%%%%%%%%%%%%%%%%%%%%%%%%%%%%%%ù
\section{Main Results}
\label{sec3}

\hh In this section, we present the main results of this paper. The asymptotic distribution of ELR test statistic under hypothesis $H_0$ will allow to build the asymptotic confidence region for the difference of the parameters of the two phases of model.  We can also test if the  models changes after observation $k$. In comparison to the obtained results for fixed $p$ (see \cite{Liu:Zou:Zhang:08}, \cite{Zi:Zou:Liu:12} for linear model,  \cite{Ciuperca:Salloum:15} for nonlinear model) where the asymptotic law is the $\chi^2$ distribution with $p$ degrees of freedom, in the case presented here, the test statistic is different and it has a standard normal asymptotic distribution. \\
In order to find this asymptotic distribution, we first need some intermediate results for studying the asymptotic behaviour of the ELR statistic.\\
We emphasize  that the  presence of the break point $k$ complicates the study and leads to a different approach in respect to a model without change-point. \\

Note that under the hypothesis $H_0$, we have: $\ez_i(\ebo)=\ez_i^0 = \XX_i \varepsilon_i$, while under $H_1$, the vector $\ez_i(\ebo)$, for $i=k+1, \cdots, n$ becomes
\begin{equation}
\label{Zib0}
\ez_i(\ebo)= \XX_i \XX_i^t (\eb^0_2-\ebo)- \ez_i^0 .
\end{equation}
When $H_0$ is true, we denote by $\es^0_n$ the matrix $\es_n(\ebo)$:
\begin{equation}
\label{S0}
\es^0_n \equiv \frac{1}{n \theta ^2 } \sum  _ {i=1} ^ {k}  \ez_i^0 \ez_i^{0t} + \frac{1}{n (1-\theta)^2} \sum _{j =k+1} ^ n \ez_j^0 \ez_j^{0t}
\end{equation}
and  by $\epsi^0_n$ the vector $\epsi_n(\ebo)$:
\begin{equation}
\label{psi}
\epsi^0_n \equiv  \frac{1}{n \theta} \sum  _ {i=1} ^ {k}  \ez_i^0 - \frac{1}{n (1-\theta )} \sum _{j =k+1} ^ n \ez_j^0.
\end{equation}

The Lagrange multiplier $\el$ is a key element in any empirical likelihood formulation. The first result concerns the convergence rate to zero of  $\el$ defined in (\ref{eq120}).   When $p$ is fixed, \cite{Zi:Zou:Liu:12} showed that $\| \el\|=O_{\eP}(n^{-1/2})$. 
When $p$ is growing along with $n$, the above rate for $\|\el\|$ is no longer valid as shown by the following proposition. In the proof we use Lemma \ref{corolary2}, Lemma \ref{lemma3} and Lemma \ref{lemma4}.
\begin{proposition}
\label{proposition1}
Suppose that assumptions (A1), (A3)-(A6) are satisfied. Then, under hypothesis $H_0$, the Lagrange multiplier $\el$ satisfies
$
\|\el\|=  O_{\eP}( p^{1/2} n^{-1/2})$. 
\end{proposition}

%%%%%%%%%%%%%%%%%%%%%%%%%%%%%%%%%%%%%%%%%%%%%%%%%%%%%%%%%%%%%%%%
%%%%%%%%%%%%%%%%%%%%%%%%%%%%%%%%%%%%%%%%%%%%%%%%%%%ù
%%%%%%%%%%%%%%%%%%%%%%%%%%%%%%%%%%%%%%%%%%%%%%%%%%%%%%%%%%%%%%%%

Accordingly to this Proposition, by assumption (A6), we have that $\|\el\| \overset{\eP} {{\longrightarrow}} 0$, as $n \rightarrow \infty$.  
More precisely, the Lagrange multiplier $\el$ has the following approximate form given by  Proposition \ref{corolary4}. 
%%%%%%%%%%%%%%%%%%%%%%%%%%%%%%%%%%%%%%%%%%%%%%%%%%%%%%%%%%%%%%%%%%%
%%%%%%%%%%%%%%%%%%%%%%%%%%%%   Corollary 4  %%%%%%%%%%%%%%%%%%%%%%%%%%%
%%%%%%%%%%%%%%%%%%%%%%%%%%%%%%%%%%%%%%%%%%%%%%%%%%%%%%%%%%%%%%%%%%%
The proof, given in Appendix, is obtained by combining Lemma \ref{lemma4},  Lemma \ref{lemma5} and Lemma \ref{lemma6}. The $p$-square matrix $\ev_n^0$ is defined by (\ref{eq4})  and the p-vector $\epsi _n^0$ by (\ref{psi}). 
\begin{proposition}
\label{corolary4}
If assumptions (A1), (A3)-(A6) are satisfied, then,  under the null hypothesis $H_0$, we have 
$ \el=(\ev_n^0)^{-1} \epsi_n^0(1+ o_{\eP}(1))$.
\end{proposition}

We prove now the following two propositions, all satisfied under hypothesis $H_0$.  
They  give two approximations for the ELR statistic  $\EL_{nk}(\ebo)$, defined by (\ref{eq6}), approximations which will allow to find its asymptotic distribution.  \\
In the proof of the following Proposition are used Lemma \ref{debut_Prop2}, Lemma \ref{lemma5}, Proposition \ref{corolary4} and Lemma \ref{lemma6}.
%%%%%%%%%%%%%%%%%%%%%%%%%%%%%%%%%%%%%%%%%%%%%%%%%%%%%%%%%%%%%%%%%%%
%%%%%%%%%%%%%%%%%%%%%%%%%   Proposition 2  %%%%%%%%%%%%%%%%%%%%%%%%
%%%%%%%%%%%%%%%%%%%%%%%%%%%%%%%%%%%%%%%%%%%%%%%%%%%%%%%%%%%%%%%%%%%
\begin{proposition}
\label{proposition2}
Suppose that assumptions (A1)-(A8) are satisfied. Then, under the null hypothesis $H_0$, we have
\begin{equation*} 
\EL_{nk}(\ebo)= n\epsi _n^{0t} (\es_n^0)^{-1} \epsi_n^0 +o_{\eP}(p^{1/2}).
\end{equation*}
\end{proposition}
For the proof of  Proposition \ref{proposition3}, given in Appendix, we use Lemma \ref{lemma4}, Lemma \ref{lemma7}, Proposition \ref{corolary4} and Proposition \ref{proposition2}.
%%%%%%%%%%%%%%%%%%%%%%%%%%%%%%%%%%%%%%%%%%%%%%%%%%%%%%%%%%%%%%%%%%%
%%%%%%%%%%%%%%%%%%%%%%%%%%%%%%%%%%%%%%%%%%%%%%%%%%%%%%%%%%%%%%%%%%%
%%%%%%%%%%%%%%%%%%%%%%%%%   Proposition 3  %%%%%%%%%%%%%%%%%%%%%%%%
%%%%%%%%%%%%%%%%%%%%%%%%%%%%%%%%%%%%%%%%%%%%%%%%%%%%%%%%%%%%%%%%%%%
%%%%%%%%%%%%%%%%%%%%%%%%%%%%%%%%%%%%%%%%%%%%%%%%%%%%%%%%%%%%%%%%%%%
\begin{proposition}
\label{proposition3}
Suppose that assumptions (A1), (A3), (A4), (A6) and (A7) are fulfilled. If the hypothesis $H_0$ is true, then we have
\begin{equation*}
\label{eq32}
\EL_n(\ebo)= n \epsi_n^{0t} (\ev^0_n)^{-1} \epsi_n^0 +o_{\eP}(p^{1/2}).
\end{equation*}

\end{proposition}

 %%%%%%%%%%%%%%%%%%%%%%%%%%%%%%%%%%%%%%%%%%%%%%%%%%%%%%%%%%%%%%%%%%%
%%%%%%%%%%%%%%%%%%%%%%%%%   Proposition 3  %%%%%%%%%%%%%%%%%%%%%%%%
%%%%%%%%%%%%%%%%%%%%%%%%%%%%%%%%%%%%%%%%%%%%%%%%%%%%%%%%%%%%%%%%%%%
%%%%%%%%%%%%%%%%%%%%%%%%%%%%%%%%%%%%%%%%%%%%%%%%%%%%%%%%%%%%%%%%%%%
The following theorem establishes the asymptotic normality of the ELR test statistic, when dimension $p$ of the explanatory variables increases to infinity as $n \rightarrow \infty$. Its proof,    given in Appendix, is very technical and moreover  the change-point presence in the model occurs in an essential way. Proposition \ref{proposition2} and Proposition \ref{proposition3} are used in the proof. 
We note that the variance of standardization $\Delta^2_n$ depends  localisation of the change in the interval $[1:n]$.

\begin{theorem}
\label{theorem2}
Under null hypothesis $H_0$, if assumptions (A1)-(A9) are satisfied and $p=o(n^{1/3})$, then \\
(i)
\begin{equation}
\label{eqt2}
\frac{n\epsi_{n} ^{0t} (\ev_n^0)^{-1} \epsi_n^0  -p}{\Delta_n /n} \overset{{\cal L}} {\underset{n \rightarrow \infty}{\longrightarrow}} {\cal N}(0,1),
\end{equation}
(ii) \hh  $\displaystyle{
\frac{\EL_{nk}(\ebo)-p}{\Delta_n /n}\overset{{\cal L}} {\underset{n \rightarrow \infty}{\longrightarrow}}{\cal N}(0,1),
}$\\
where $\Delta^2_n= \sum  ^{n}_{i=1} \sigma_i^2$, with $\sigma_1^2 =\theta^{-4} \big( ( \XX_1 ^t (\ev_n^0)^{-1} \XX_1 )^2 \eE[\varepsilon_1^4]-[\tr ((\ev_n^0)^{-1} \ev^0_{(1)}) ]^2\big)$ and:\\
- for  $i=2, \cdots, k+1$,
\[
 \sigma^2_i =\frac{4}{\theta^4} \sum  ^{i-1}_{l=1} \tr \big((\ev_n^0)^{-1}  \ev^0_{(i)}  (\ev_n^0)^{-1} \ev^0_{(l)} \big) 
+\frac{( \XX_i ^t (\ev_n^0)^{-1} \XX_i )^2 \eE[\varepsilon_1^4] 
-   [\tr ((\ev_n^0)^{-1} \ev^0_{(i)}) ]^2}{\theta^4},
\]
- for  $ i=k+2, \cdots ,  n$,
\begin{eqnarray*}
 \sigma^2_i &=& \frac{4}{\theta^2(1-\theta)^2} \sum  ^{k}_{l=1}  \tr \big((\ev_n^0)^{-1}  \ev^0_{(i)}  (\ev_n^0)^{-1} \ev^0_{(l)} \big) 
+\frac{4}{(1-\theta)^4} \sum  _{l=k+1}^{i-1}   \tr \big((\ev_n^0)^{-1}  \ev^0_{(i)}  (\ev_n^0)^{-1} \ev^0_{(l)} \big)  
\nonumber \\ &&
+\frac{( \XX_i ^t (\ev_n^0)^{-1} \XX_i )^2 \eE[\varepsilon_1^4] 
-   [\tr ((\ev_n^0)^{-1} \ev^0_{(i)}) ]^2}{(1-\theta)^4}.
\end{eqnarray*}
\end{theorem}

The following result is an immediate corollary of Theorem \ref{theorem2}.

%%%%%%%%%%%%%%%%%%%%%%%%%%%%%%%%%%%%%%%%%%%%%%%%%%%%%%%%%%%%%%%%
%%%%%%%%%%%%%%%%%%%%%%%%%%%Remarque 1%%%%%%%%%%%%%%%%%%%%%%%%
%%%%%%%%%%%%%%%%%%%%%%%%%%%%%%%%%%%%%%%%%%%%%%%%%%%%%%%%%%%%%%%%
\begin{Corollary}
\label{remarque1}
Testing the null hypothesis $H_0 : \eb=\eb_2=\ebo$ against the alternative hypothesis $H_1 : \eb=\ebo, \, \eb_2 \neq \ebo$, is equivalent to constructing the confidence regions for $\del=\ebo-\eb^0_2$, or to testing the null hypothesis $H'_0 : \del=\textbf{0}_p$. Then, based to Theorem \ref{theorem2}, in order to test $H_0$ against $H_1$, we consider the following asymptotic test statistic 
\begin{equation}
\label{Zb}
{\cal Z}(\ebo) \equiv \frac{n\epsi _n^{t}(\ebo) \big(\ev_n^0\big)^{-1} \epsi_n(\ebo)-p}{\Delta_n /n}.
 \end{equation}
\end{Corollary}

Note that $\epsi_n(\ebo) $ through $\ez_i(\ebo)$, given by relation (\ref{Zib0}), for $i= k+1, \cdots, n$, depends of $\del=\ebo-\eb^0_2$.\\

The asymptotic behaviour  under hypothesis $H_1$  of the test statistic ${\cal Z}(\ebo)$ is given by the following theorem. We show that ${\cal Z}(\ebo)$ diverges under alternative hypothesis.
\begin{theorem}
\label{theorem_H1}
Under alternative hypothesis $H_1$, if assumptions (A1)-(A9) are satisfied and $p=o(n^{1/3})$, then $|{\cal Z}(\ebo) |\overset{\eP} {\underset{n \rightarrow \infty}{\longrightarrow}} \infty$.
\end{theorem}

Theorem \ref{theorem2} and Theorem \ref{theorem_H1}  allow to build the asymptotic confidence region for the parameter $ \del = \ebo- \eb_2$.

\begin{Corollary}
\label{remarque1b}
 The $\alpha$-level asymptotic confidence region for $\del$ is 
\begin{equation}
\label{intr}
{\cal R}_{1- {\alpha}/{2}}= \left \{ \del :   \big |{\cal Z}(\ebo) \big  | < c_{1-{\alpha}/{2}} \right \},
\end{equation}
where $c_{1- {\alpha}/{2}}$ is the quantile of the standard normal distribution. 
\end{Corollary}

For simulations, in order to calculate ${\cal R}_{1- {\alpha}/{2}}$, the matrix $\ev_n^0$ is  firstly calculated  by relation (\ref{eq4}). Once the model has been generated, we calculate $\ez_i^0$ and then $\ez_i(\ebo)$ by relation (\ref{Zib0}). Finally, we calculate $\epsi_n(\ebo)$ by (\ref{psib}), $\Delta_n$ by Theorem \ref{theorem2} and the test statistic ${\cal Z}(\ebo)$ by (\ref{Zb}). For  $M$ Monte Carlo replications of the model, the coverage rate (CR),  is  the number of times when $| {\cal Z}(\ebo)|$ is less than $c_{1-{\alpha}/{2}}$, divided by $M$.\\
For applications on real data, we will test model (\ref{eq1}) against model (\ref{eq2}). For these models, we know $n$ values for the response variable $Y$ and for the $p-1$ explanatory variables $X_2, \cdots , X_p$. The point $k$, where we want to test if there is a change, is known, while the values of $\ebo$ on the first phase can be unknown. \\
If $\ebo$ is unknown, then it is estimated by  a convergent estimator on the observations $i=1, \cdots, k$, for example by LS method or quantile method, depending on the   distribution of $Y$. Once we dispose of an estimator $\hat \eb_k$ for $\ebo$, the variance $\sigma ^2$ of $\varepsilon$ is estimated afterwards by a convergent estimator,  for instance $\hat \sigma^2_k = (k-p)^{-1} \sum^k_{i=1}(Y_i - \XX^t_i \hat  \eb_k)^2$. We calculate thereby $\ev_{(i)}^0= \hat \sigma^2_k \XX_i \XX_i^t$, for any $i=1, \cdots, n$ and then $\ev_n^0$ by (\ref{eq4}). For any $i=1, \cdots, n$ we calculate $\ez_i(\hat \eb_k)= \XX_i(Y_i -\XX_i^t \hat \eb_k)$, which will allow us to calculate the vector  $\epsi_n(\hat \eb_k)$ of  relation  (\ref{psib}). With all of these elements in place, we can calculate the value of the statistic ${\cal Z}(\hat \eb_k) = \big({n\epsi _n^{t}(\hat \eb_k) \big(\ev_n^0\big)^{-1} \epsi_n(\hat \eb_k)-p}\big)\left(\Delta_n /n\right)^{-1}$,  using for $\Delta_n$ the  relation given in  Theorem \ref{theorem2}. For a given size $\alpha \in (0,1)$, if the value of $| {\cal Z}(\hat \eb_k)|$ is less than $c_{1-{\alpha}/{2}}$, then hypothesis $H_0$ is accepted, that is to say that the model does not change after observation $k$, otherwise hypothesis $H_1$ is accepted.\\
If $\ebo$ is known, we can consider as an estimator for $\sigma^2$: $\hat \sigma^2_k = (k-p)^{-1} \sum^k_{i=1}(Y_i - \XX^t_i   \ebo)^2$. For any $i=1, \cdots, n$ we calculate $\ez_i(\ebo)= \XX_i(Y_i -\XX_i^t  \ebo)$ and afterwards $\epsi_n( \ebo)$ by relation  (\ref{psib}). Finally, the absolute value of ${\cal Z}(\ebo)$ will be compared with $c_{1-{\alpha}/{2}}$.

\begin{Remark}
Compared to  \cite{Liu:Zou:Zhang:08}, where, for fixed $p$, a test statistic is proposed for testing the presence of the change-point, by maximizing ELR in respect to $\eb$ and $\el$, in the present work we fix the parameter on the first phase and we test  whether the parameter of the second phase is the same. In  \cite{Liu:Zou:Zhang:08}, the system of equations in $\el$ and $\eb$ of the score functions must be solved, which can be numerically quite tedious. In this paper, apart from the fact that we consider $p \rightarrow \infty$, using theoretical properties for the Lagrange multiplier $\el$, we propose a simpler form for ELR statistic, easier to use in practice. Parameter $\ebo$, if it is unknown, can be estimated on the observations $1, \cdots, k$ by a simpler computational method, in order to obtain $\hat \eb_k$ a convergent estimator, i.e. $\| \hat \eb_k - \ebo\| =o_{\eP}(1)$.
% As we shall see in the next section, where simulations are presented, if we consider a similar statistic, as (\ref{ref17}), but without taking into account that $p \rightarrow \infty$, then it follows asymptotically a chi-square distribution with $p$ degrees of freedom. Statistic (\ref{ref17}) will give worse results than statistic (\ref{Zb}) in terms of coverage accuracy.  
\end{Remark}

%%%%%%%%%%%%%%%%%%%%%%%%%%%%%%%%%%%%%%%%%%%%%%%%%%%%%%%%%%%%%%%%%%%
%%%%%%%%%%%%%%%%%%%%%%%%%%%%%%%%%%%%%%%%%%%%%%%%%%%%%%%%%%%%%%%%%%%
%%%%%%%%%%%%%%%%%%%%%%%%%%%%%%%%%%%%%%%%%%%%%%%%%%%%%%%%%%%%%%%%%%%
     %%%%%%%%%%%%%%%%%%%%   Simulations  %%%%%%%%%%%%%%%%%%%%%%%%
%%%%%%%%%%%%%%%%%%%%%%%%%%%%%%%%%%%%%%%%%%%%%%%%%%%%%%%%%%%%%%%%%%%
%%%%%%%%%%%%%%%%%%%%%%%%%%%%%%%%%%%%%%%%%%%%%%%%%%%%%%%%%%%%%%%%%%%
%%%%%%%%%%%%%%%%%%%%%%%%%%%%%%%%%%%%%%%%%%%%%%%%%%%%%%%%%%%%%%%%%%%
\section{Simulation study}
\label{sec4}
\hh We now conduct simulation studies to evaluate, in terms of coverage accuracy and empirical power, the test statistic specified by Theorem \ref{theorem2}(ii), with $\EL_{nk}(\ebo)$ approximated by Proposition \ref{proposition3}, ie ${\cal Z}(\ebo)$ given by relation (\ref{Zb}).  
For these studies, we use Monte Carlo simulations. Throughout, we consider the  size $\alpha=0.05$. \\
\hh The $p$ explanatory variables are generated as follows: $X_1=1$ and $(X_2, \cdots, X_p ) \sim {\cal N}_{p-1}(\textbf{0}_{p-1}, \bf{\Sigma})$, with the covariance matrix ${\bf \Sigma}=(\varsigma_{hl})$, $\varsigma_{hl}=2^{-|h-l|}$, $1 \leq h, l \leq p-1$, the same matrix considered by \cite{Guo:Zou:Wang:Chen:13},  for a model without change-point. In order to be in a fixed design, we consider the same realization for  $(\XX_i)_{1 \leq i \leq n}$ for each Monte Carlo replication. \\
Concerning the coefficients $\ebo$ of the model, under $H_0$, we take  $\ebo=(\beta^0_1,\beta^0_2,  \cdots , \beta^0_{p})=(1, 2 , \cdots, p)$. \\
 For  model (\ref{eq2}), under hypothesis $H_0: \eb =\eb_2 =\eb^0$, we first calculate the  coverage rate (CR) based on Corollary \ref{remarque1b}, for a given change-point $k$.\\
We consider different values for $n$ and $k$ and two different distributions for the errors $(\varepsilon_i)$: standard normal distribution  ${\cal N}(0,1)$ and    $\varepsilon \sim {\cal E}xp(1)-1$, where ${\cal E}xp(1)$ is the exponential distribution with mean $1$. 
\subsection{Importance of assumptions (A5), (A6)}
\hh In this subsection we realise throughout 2000 Monte Carlo replications for studying the behaviour of the test statistic behaviour, under null hypothesis and afterwards, when model has a change-point. Coverage rate and empirical power are investigated. Values of $n$ and $k$ are   $n \in \{ 20, 100, 200, 400, 600, 1000\}$,  $k \in \{5, 25, 75, 280, 350 \}$.

%%%%%%%%%%%%%%%%%%%%%%%%%%%%%%%%%%%%%%%%%%ù
\begin{table}[h!]
\begin{center}
\caption{Coverage rates (CR)   for Exponential and Gaussian errors.}
\label{tabl1}
\ \\
\begin{tabular}{|c c c ||c||c| }   \hline
& & &  \multicolumn{1}{c||}{ Exponential errors } &  \multicolumn{1}{c|}{ Gaussian errors } \\
    $n$  & $k$ &  $p$  &   CR &  CR   \\ \noalign{\smallskip}\hline\noalign{\smallskip}
  20& 5& 2 &0.94  & 0.96 \\
     & & 3 &0.92  & 0.95 \\
     & &  5 &0.91  & 0.92\\
      & &  7 &0.77  & 0.87\\
   \hline
    100& 25& 2 &0.94  & 0.95 \\
     & & 10 &0.88  & 0.94 \\
      & & 20 & 0.71  & 0.84 \\
    \hline
    200& 75& 2 &0.93  & 0.94\\
        & & 10 &0.89  & 0.93\\
     &  & 20 & 0.77 &0.94\\
      & & 50 & 0.73 &  0.85 \\
         \hline
  %%%%%%%%%%%%%%%%%%%%%%%%%%%%%%%%%%%%%%%%%%%%%%%%%
  %%%%%%%%%%%%%%%%%%%%%%%%%%%%%%%%%%%%%%%%%%%%%%%%%%   
      400& 280& 2 &0.92  & 0.91 \\
      & & 10 &0.87  & 0.88 \\
     & &  20 & 0.79 &  0.86\\
      & & 50 &0.75 &  0.84 \\
       & &  100&  0.60  & 0.83\\
       \hline
  %%%%%%%%%%%%%%%%%%%%%%%%%%%%%%%%%%%%%%%%%%%%%%%%%
  %%%%%%%%%%%%%%%%%%%%%%%%%%%%%%%%%%%%%%%%%%%%%%%%%%  
    600& 350& 2 & 0.94  & 0.92 \\
      & & 10 & 0.92  & 0.93 \\
     & &  20 &  0.86  & 0.87 \\
      & &  50 & 0.85 & 0.88 \\
       & &  100  & 0.78&  0.84 \\
       & &  200 &0.66  & 0.84 \\
        & &  300 &0.54  & 0.70           \\ \hline
  %%%%%%%%%%%%%%%%%%%%%%%%%%%%%%%%%%%%%%%%%%%%%%%%%
  1000& 350& 2 &0.95 &  0.94   \\
  &  & 50 &0.86  &   0.92  \\
  &  & 100 &0.84  &   0.89   \\
   &  & 200 &0.78 &   0.80  \\
   &  & 300 &0.52 &   0.77  \\
  %%%%%%%%%%%%%%%%%%%%%%%%%%%%%%%%%%%%%%%%%%%%%%%%%%  
    \hline
\end{tabular}
\end{center}
\end{table}
\textit{Analyse of coverage rate.} The results are summarized in Table \ref{tabl1} where we give  \textit{CR = 1 - empirical size},  based on Corollary \ref{remarque1b}, relation (\ref{Zb}). For $n$ and $k$ fixed, the CRs decrease when $p$ increases, this decreasing trend being more pronounced  in the exponential error case. 
We observe that whether for exponential errors or for gaussian errors, if  assumptions (A5), (A6) are not satisfied, then the CRs are well below 0.95. 
 These results are in accordance with those obtained by \cite{Guo:Zou:Wang:Chen:13}, for models without change-point, with fixed design. \\
In order to confirm this supposition, in Table \ref{tabl2}, the values of $n$ and $k$ are varied such that $\theta=k/n =1/2$ and $p$ satisfies (A5), (A6). We obtain then that the CRs are larger than 0.90. \\
\hh \textit{Analyse of power.} Under $H_1$, we consider  $\eb_2^0=1-\ebo$. In the all considered cases, for $n$, $k$, $p$ and $\varepsilon$ in Tables \ref{tabl1} and \ref{tabl2},  we obtain that the empirical powers are equal to 1. 

\begin{table}[h!]
\begin{center}
\caption{Coverage rates (CR),  by 2000 Monte Carlo replications,  for Exponential and Gaussian errors,  $\theta=1/2$.}
\label{tabl2}
\ \\
\begin{tabular}{|c c c ||c||c| }   \hline
& & &  \multicolumn{1}{c||}{ Exponential errors } &  \multicolumn{1}{c|}{ Gaussian errors } \\
    $n$  & $k$ &  $p$  &   CR &  CR   \\ \noalign{\smallskip}\hline\noalign{\smallskip}
  20& 10& 2 &0.92  & 0.96 \\
   100  & 50 & 4 &0.93  & 0.94 \\
    200 & 100 &  5 &0.96  & 0.97\\
    400  & 200 & 10 &0.90  & 0.94 \\
       600& 300 & 20 &0.90  & 0.92 \\
  800   & 400 & 20 &0.93  & 0.93 \\
    800  & 400& 30 & 0.90  & 0.91 \\
    2000& 1000& 30 &0.91  & 0.91\\
       \hline
\end{tabular}
\end{center}
\end{table} 
\subsection{CR's improvement}
\hh In order to obtain more precise false probabilities, for fixed size $\alpha$, we will calculate, by 10000 Monte Carlo replications, the $(1 -\alpha/2)$ and $\alpha /2 $ quantiles,  denoted $\hat {c}_{1}$, $\hat {c}_{2}$, respectively, for  test statistic ${\cal Z}(\ebo)$.  We consider the new critical value $\hat c_{1- \alpha / 2} \equiv \max (\hat {c}_{1}, |\hat {c}_{2}|)$.  These new critical values, for $p=50$, for different values of $n$ and $k$,  for ${\cal N}(0,1)$ and ${\cal E}xp(1)-1$ distribution  errors,  are given in Table \ref{tabl3}. These values are not influenced by value of $k$, for fixed  $n$. This is observed by calculating  $\hat c_{1- \alpha / 2}$ for $k$ such that $\theta_1=k/n=3/8$  and afterwards we calculate the CRs, denoted $\widehat {CR} $, for another $k$ such that $\theta_2=k/n=5/8$.  We observe that the values of $\hat c_{1- \alpha / 2}$ are larger than the   quantile of the standard normal distribution and $\hat c_{1- \alpha / 2}$ are larger for exponential errors than those for normal errors. On the other hand, the values of $\hat c_{1- \alpha / 2}$ decrease when $n$ (and $k$) increases and they approach to quantile of ${\cal N}(0,1)$.\\
  In the same Table, are given empirical powers, denoted $\hat \pi$, calculated for $\eb_2^0=1-\ebo$ under $H_1$, considering $\hat c_{1- \alpha / 2}$ as critical value. We obtain that all $\hat \pi$ are equal to 1. \\
  If under $H_1$, only two components of $\ebo$ change: $\eb_2^0$ is such that $\beta^0_{2,j}=\beta^0_{j}$ for all $j \in \{ 1, \cdots, p\} \setminus \{ 3, 30\}$, $\beta^0_{2,3}=\beta^0_{3}+1$, $\beta^0_{2,30}=\beta^0_{30}+1$, we always get $\hat \pi=1$. Here we have denoted by   $\beta^0_{2,j}$ the $j$th component of $\eb^0_2$. Hence, even if there is a small change in the coefficients, most coefficients remaining unchanged, the test statistic detects this change. 
\begin{table}[h!]
\begin{center}
\caption{Empirical critical value $\hat c_{1- \alpha/2}$ and corresponding   coverage rates ($\widehat{CR}$), empirical powers ($\hat \pi$),  for Exponential and Gaussian errors, $p=50$, $\eb_2^0=1 -\ebo$.}
\label{tabl3}
\ \\
\begin{tabular}{|cc ||c|cc||c|cc| }   \hline
&  &   \multicolumn{3}{c||}{ Exponential errors } &  \multicolumn{3}{c|}{ Gaussian errors } \\
    $n$   &  $k$   & $\hat c_{1- \alpha /2}$ & $\widehat{CR}$ &  $\hat \pi$  & $\hat c_{1- \alpha}$ & $\widehat{CR}$ &  $\hat \pi$ \\ \noalign{\smallskip}\hline\noalign{\smallskip}
     %%%% 14 avril, 10000 MC n=200
  200  & 75 &  4.01 & 0.97 & 1& 3.24 & 0.97 & 1 \\
     &     125 &     & 0.97 & 1&   & 0.93 & 1 \\ \hline 
     %%%% 14 avril, 10000 MC n=400
      400  &  150 &  3.41  & 0.97  & 1 &  2.96 & 0.97  & 1  \\
     &     250 &      & 0.97  & 1 &   &  0.98 & 1  \\ \hline 
     %%%% 13 avril: 10000MC n=600 
      600 &  225 &  3.40  & 0.97  & 1 & 2.85  &  0.97 & 1 \\
     &     375 &      & 0.97  & 1 &   &  0.98 & 1 \\ \hline 
%%%%%% 13(exp) et 17 (N) avril: simus pour 10000 MC n=800
         800 &    300 &   2.68  & 0.95  & 1 & 2.43  & 0.97  &  1 \\
     &      500 &       & 0.95  & 1 &   & 0.95  & 1 \\ \hline 
     %%%%%% 13 avril: simus pour 4000 MC n=2000
       2000  &  750 &   2.48   & 0.97  & 1 &  2.30 &  0.97 & 1 \\
            &    1250 &       & 0.95  & 1 &   & 0.93  & 1 \\ \hline 
\end{tabular}
\end{center}
\end{table} 
\subsection{Conclusion of simulations}
\hh Proposed test statistic (\ref{Zb}), with ${\cal N}(0,1)$ the asymptotic distribution under $H_0$, involves the construction of a confidence region for the parameters of the  second phase of the model (on  observations $k+1, \cdots, n$). \\
If assumptions (A5), (A6) are satisfied, then the coverage rates are close to the nominal coverage level. Contrariwise, if the coefficients change on the second phase, the test always detects this change. For improving the coverage rate in the case $n-k  \not \gg p^2$ or $k \not \gg p^2$, we proposed to calculate new critical values.  With these critical values, the rate of false changes is generally smaller than the size $\alpha$. If there are changes in the coefficients of the second phase of the model, the test statistic based on the new confidence region always  detects this change. For fixed $p$, if $n$ and $k$ increase, such that $k/n$=constant, then these new critical values decrease and approach  the $(1- \alpha)$ quantile of ${\cal N}(0,1)$ distribution. 
%%%%%%%%%%%%%%%%%%%%%%%%%%%%%%%%%%%%%%%%%%%%%%%%%%%%%%%%%%%%%%%%%%%
%%%%%%%%%%%%%%%%%%%%%%%%%%%%%%%%%%%%%%%%%%%%%%%%%%%%%%%%%%%%%%%%%%%
%%%%%%%%%%%%%%%%%%%%%%%%%%%%%%%%%%%%%%%%%%%%%%%%%%%%%%%%%%%%%%%%%%%
     %%%%%%%%%%%%%%%%%%%%   Appendix  %%%%%%%%%%%%%%%%%%%%%%%%
%%%%%%%%%%%%%%%%%%%%%%%%%%%%%%%%%%%%%%%%%%%%%%%%%%%%%%%%%%%%%%%%%%%
%%%%%%%%%%%%%%%%%%%%%%%%%%%%%%%%%%%%%%%%%%%%%%%%%%%%%%%%%%%%%%%%%%%
%%%%%%%%%%%%%%%%%%%%%%%%%%%%%%%%%%%%%%%%%%%%%%%%%%%%%%%%%%%%%%%%%%%
\section{Appendix}
\label{sec5}
%%%%%%%%%%%%%%%%%%%%%%%%%%%%%%%%%%%%%
%%%%%%%%%%%%%%%%%%%%%%%%%%%%%%%%%%%%%
%%%%%%%%%%%%%%%%%%%%%%%%%%%%%%%%%%%%%
\hh This section is divided into two subsections. In the first we give the proofs of the Propositions and of the Theorems. In the second subsection, we present Lemmas ans their proofs. \\

We recall that under the hypothesis $H_0$, the vector   $\ez_i(\ebo) $ is $\ez_i^0=\XX_i \varepsilon_i$. Then, in the all proofs, if hypothesis  $H_0$ is true, we will use $\ez_i^0$ instead of $\ez_i(\ebo) $.

\subsection{Proposition and Theorem proofs}
 
\noindent {\bf Proof of Proposition \ref{proposition1}.}  Let us write $\el$ as $\el = \|\el\| \textbf{u}$, where $\textbf{u}$ is a $p$-vector with norm one. Using Lemma \ref{lemma3}, for $\theta=k/n$,   we have, for any $i=1,\ldots, k$, with probability one
\begin{equation}
\label{maj1}
0< 1+ \frac{\|\el\|}{\theta} \textbf{u}^t \ez_i^0 \leq 1+ \|\el\| T_n^0 ,
\end{equation}
where $ T_n^0 \equiv \max _{\substack{1\leqslant i \leqslant k\\  k+1 \leqslant j \leqslant n}} \big \{ \displaystyle (k/n)^{-1}\|\ez_i ^0\|, (1-k/n)^{-1}\|\ez_j ^0\|  \big \}$. 
For $j=k+1,\ldots,n $ we have, with probability 1,
\begin{equation}
\label{maj2}
0< 1- \frac{\|\el\|}{1-\theta} \textbf{u}^t \ez_j ^0 \leq 1+ \|\el\|T_n^0.
\end{equation}
Using relations (\ref{maj1}) and (\ref{maj2}), then we get from (\ref{eq120}) that 
\begin{eqnarray*}
\label{eq15}
0&=& \frac{1}{n} \sum  _ {i=1} ^ {k}  \frac{\textbf{u}^t  \ez_i^0}{\theta +\el^t \ez_i^0} - \frac{1}{n} \sum  _{j =k+1} ^ n \frac{\textbf{u}^t  \ez_j^0}{1- \theta -\el^t \ez_j^0} 
\nonumber \\ 
&=&   \frac{1}{n\theta} \sum  _{i=1} ^{k} \textbf{u}^t \ez_i^0 - \frac{1}{n\theta^2} \| \el\| \sum  _{i=1} ^{k} \frac{\textbf{u}^t \ez_i^0  \ez_i^{0t} \textbf{u }}{1+\theta^{-1}{\|\el\|}\textbf{u}^t \ez_i^0  }
%\nonumber \\ &&
-\frac{1}{n(1-\theta)} \sum  _{j=k+1} ^{n} \textbf{u}^t \ez_j^0 - \frac{1}{n(1-\theta)^2} \| \el\| \sum  _{j=k+1} ^{n} \frac{\textbf{u}^t \ez_j^0  \ez_j^{0t} \textbf{u} }{1- (1-\theta)^{-1}{\|\el\|}\textbf{u}^t \ez_j^0  }. 
\end{eqnarray*}
By the last equality, using also notations given by (\ref{S0}) and (\ref{psi}), it follows that 
\[0 \leq \textbf{u}^t \epsi_n^0 - \frac{\|\el\| }{1+\|\el\|T_n ^0 } \textbf{ u}^t \es_n^0 \textbf{u} .\]
 Then, we have  with probability one, that $\textbf{u}^t \epsi_n ^0 (1+\|\el\|T_n^0  ) \geq \|\el\| \textbf{u}^t \es_n^0 \textbf{u}$. Therefore 
\begin{equation}
\label{eq16}
\|\el\| \leq \frac{\textbf{u}^t \epsi_n^0}{ \textbf{u}^t \es_n^0 \textbf{u}-\textbf{u}^t \epsi_n^0 T_n^0}.
\end{equation}
On the other hand, we have $|\textbf{u}^t \epsi_n^0| \leq \|\epsi_n^0 \|$. Then, using  Lemma \ref{lemma4}, we obtain that $\|\epsi_n^0 \|=O_{\eP}(n^{-1/2} p^{1/2})$, which gives 
\begin{equation}
\label{eq17}
\textbf{u}^t \epsi_n ^0=O_{\eP}(n^{-1/2} p^{1/2}).
\end{equation}
Using Lemma \ref{lemma3} and relation (\ref{eq17}), we have that $\textbf{u}^t \epsi_n^0 T_n^0= O_{\eP}(n^{-1/2} p^{1/2}) o_{\eP}(n^{1/q} p^{1/2})= o_{\eP}(n^{(-q+2)/2q} p)$. Then, by assumption (A5), we obtain that 
$
\textbf{u}^t \epsi_n^0 T_n^0=o_{\eP}(1)$.\\ 
On the other hand, according to Lemma \ref{corolary2}, $\textbf{u}^t \es_n^0 \textbf{u } \geq \gamma_p(\es_n^0) > C_0>0$ holds with a probability tending to 1 as $n \rightarrow \infty$. Then, for  relation (\ref{eq16}), we obtain that 
\[ \|\el\|= O_{\eP}( |\textbf{u}^t \epsi_n^0|/C_0)=O_{\eP}(\| \epsi_n^0 \|)= O_{\eP}(p^{1/2}n^{-1/2} ).\]
\hspace*{\fill}$\blacksquare$ \\

%%%%%%%%%%%%%%%%%%%%%%%%%%%%%%%%%%%%%%%%%%%%%%%%%%%%%%%%%%%%%%%%
\noindent\textbf{Proof of Proposition \ref{corolary4}.} By Lemma \ref{lemma5}  we have  that $\el= (\es_n^0)^{-1} \left(\er_n^0+ \epsi_n^0\right)(1+ o_{\eP}(1))$. \\
In the other hand, by Lemma \ref{lemma4}, we have that $\|\epsi_n^0 \|=O_{\eP}( p^{1/2 }n^{-1/2} )$. Using this fact and relation (\ref{eq23}), we obtain %relation (\ref{eq25}) becomes
\begin{eqnarray*}
\label{eq26}
\frac{ \|\er_n ^0 \|}{\|\epsi_n ^0 \|}= o_{\eP}(n^{1/q} p^{1/2} n^{-1/2} p^{1/2} ) 
= o_{\eP}(pn^{(2-q)/2q} ).
\end{eqnarray*}
Therefore, by assumption (A5) we obtain that $\er_n^0 =\epsi_n^0  o_{\eP}(1).$ 
Then $\el=(\es_n^0)^{-1} \epsi_n^0 (1+o_{\eP}(1)). $ \\ 
In the other hand, by Lemma \ref{lemma6}, we have that $\left((\es_n^0)^{-1} -(\ev_n^0)^{-1}\right) \epsi_n^0=(\ev_n^0)^{-1} \epsi_n^0 o_p(1)$. Then, $(\es_n^0)^{-1}  \epsi_n^0= (\ev_n^0)^{-1} \epsi_n^0 (1+ o_{\eP}(1))$. Therefore, for $\el=(\es_n^0)^{-1} \epsi_n^0(1+o_{\eP}(1))$, we obtain that 
\begin{equation}
\label{eq31}
\el=(\ev_n^0)^{-1} \epsi_n^0 (1+ o_{\eP}(1)).
\end{equation}
\hspace*{\fill}$\blacksquare$ \\
%%%%%%%%%%%%%%%%%%%%%%%%%%%%%%%%%%%%%%%%%%%%%%%%%%%%%%%%%%%%%%%%
%%%%%%%%%%%%%%%%%%%%%%%%%%%%%%%%%%%%%%%%%%%%%%%%%%%%%%%%%%%%%%%%
 
\noindent {\bf Proof of Proposition \ref{proposition2}.} 
By Lemma \ref{debut_Prop2} we have that 
\begin{equation}
\label{etiqq}
\EL_{nk}(\ebo)= 2n \el ^t\epsi_n^0 -n\el^t \es_n^0 \el +{\cal E}_3+o_{\eP}(1), 
\end{equation}
with ${\cal E}_3 \equiv \displaystyle \frac{2}{3} \big (  \displaystyle \frac{1}{\theta^{3}}  \sum  _ {i=1} ^ {k}  (\ez_i^{0t} \el)^3-  \displaystyle \frac{1}{(1-\theta)^{3}}  \sum _{j =k+1} ^ n  (\ez_j^{0t}\el)^3\big) $. \\
Consider now, the following $p$-vector
\begin{equation}
\label{R0n}
\er_n^0  \equiv  \frac{1}{n\theta^3}  \sum  _ {i=1} ^ {k}  \ez_i^{0}  ( \el ^t \ez_i^0)^2 - \frac{1}{n(1-\theta)^3} \sum  _ {j =k+1} ^ n   \ez_j^{0}  ( \el ^t \ez_j^0)^2.
\end{equation}
By Lemma \ref{lemma5} we have that $\el =(\es^0_n)^{-1}(\er_n^0+\epsi^0_n)(1+o_{\eP}(1))$. Then, we have for (\ref{etiqq}) that
\begin{equation}
\label{eq170} 
\EL_{nk}(\ebo)=n \epsi_n ^{0t} (\es_n^0)^{-1}  \epsi_n^0 - n\er_n^{0t} (\es_n^0)^{-1} \er_n^0 +{\cal E}_3+o_{\eP}(1).
\end{equation}
We now study ${\cal E}_3$ and $n(\er_n^0)^t (\es_n^0)^{-1} \er_n^0$ in parallel. By Proposition \ref{corolary4}, we have that $\el=(\ev^0_n)^{-1}\epsi_n ^0(1+o_{\eP}(1))$, which implies
 \begin{equation}
\label{eqz}
\ez_i^{0t} \el = \ez_i^{0t}(\ev^0_n)^{-1} \epsi_n^0(1+o_{\eP}(1)).
\end{equation}
Then, ${\cal E}_3$ becomes
\begin{eqnarray}
\label{F1}
{\cal E}_3&=&\frac{2}{3\theta^3} \sum  _ {i=1} ^ {k}  \bigg(\ez_i^{0t} (\ev^0_n)^{-1/2}(\ev_n^0)^{-1/2} \big (\frac{1}{n \theta} \sum  _ {i=1} ^ {k}  \ez_i^0 + \frac{1}{n(\theta-1)} \sum _{j =k+1} ^ n  \ez_j^{0} \big ) \bigg )^3(1+o_{\eP}(1))
 \nonumber \\ &&
+ \frac{2}{3(\theta-1)^3}  \sum _{j =k+1} ^ n   \bigg( \ez_j^{0t} (\ev^0_n)^{-1/2} (\ev^0_n)^{-1/2} \big(\frac{1}{n \theta} \sum  _ {i=1} ^ {k}  \ez_i^0 + \frac{1}{n(\theta-1)} \sum _{j =k+1} ^ n  \ez_j^{0} \big ) \bigg)^3   (1+o_{\eP}(1)). 
 \nonumber 
\end{eqnarray}
Using  notations given by (\ref{eq8}), (\ref{eq9}) and the strong law of large numbers (Markov's Theorem), we obtain
\begin{eqnarray}
\label{F33}
{\cal E}_3&=&\frac{2}{\theta^3} \sum  _ {l=1} ^ {k}  \big (\ew_l^{0t} \big(\frac{1}{n \theta} \sum  _ {i=1} ^ {k}\ew_i^0+\frac{1}{n(\theta-1)} \sum _{j =k+1} ^ n 
\ew_j^{0} \big)\big )^3 
% \nonumber \\ &&
+ \frac{2}{(1-\theta)^3}  \sum _{l =k+1} ^ n   \big (\ew_l^{0t} \big(\frac{1}{n \theta} \sum  _ {i=1} ^ {k}  \ew_i^0+\frac{1}{n( \theta -1)} \sum _{j =k +1} ^ n   \ew_j^{0} \big) \big)^3 (1+o_{\eP}(1)) \nonumber \\ 
&=& \frac{2n}{3} \sum  _ {r,s,u=1} ^ {p} \omega ^r \omega ^s \omega ^u \big(\frac{1}{n\theta^3} \sum  _ {i=1} ^ {k}  w_{i,r}^0 w_{i,s}^0 w_{i,u}^0 + \frac{1}{n(\theta-1)^3}  \sum _{j =k+1} ^ n   w_{j,r}^0 w_{j,s}^0 w_{j,u}^0    \big) (1+o_{\eP}(1)) \nonumber \\ 
&=& \frac{2n}{3} \sum  _ {r,s,u=1} ^ {p} \omega ^r \omega ^s \omega ^u \alpha^{rsu} (1+o_{\eP}(1)).
\end{eqnarray}
In the other hand,   replacing $\el$ in  relation (\ref{R0n})  we obtain 
\begin{equation}
\label{rari}
\er_n^0= \bigg ( \frac{1}{n\theta^3} \sum  _ {i=1} ^ {k}  \ez_i^0  \big (\ez_i^{0t} \epsi_n ^{0t}(\ev^0_n)^{-1} \big)^2-\frac{1}{n(1-\theta)^3}  \sum  _ {j =k+1} ^ n  \ez_j^0  \big (\ez_j^{0t} \epsi_n ^{0t}(\ev^0_n)^{-1} \big)^2\bigg ) (1+o_{\eP}(1)).
\end{equation}
Consider now for   $n(\er^0_n)^t (\es^0_n)^{-1} \er^0_n$  the following decomposition
\begin{equation}
\label{eq127}
 \er^{0t}_n (\es^0_n)^{-1} \er^0_n
= \er^{0t}_n \big((\es^0_n)^{-1}-(\ev^0_n)^{-1} \big) \er_n^0 +  \er^{0t}_n (\ev^0_n)^{-1} \er^0_n.
\end{equation} 
By Lemma \ref{lemma6}(ii), we have that $\big((\es^0_n)^{-1}-(\ev^0_n)^{-1}\big) \er_n^0 = (\ev^0_n)^{-1} \er_n^0 o_{\eP}(1)$. Then,  relation (\ref{eq127}) becomes
%\begin{equation}
%\label{eq128}
$ \er^{0t}_n (\es^0_n)^{-1} \er^0_n= \er^{0t}_n (\ev^0_n)^{-1}  \er^0_n(1+o_{\eP}(1))$.
% \end{equation}
 Using  relations (\ref{psi}), (\ref{R0n}), (\ref{rari}) and the fact that $\ew_i^0= (\ev^0_n)^{-1/2} \ez_i^0$, for $ i=1, \cdots, n$, we have that $n \er^{0t}_n (\es^0_n)^{-1} \er^0_n$  can be  written
\begin{eqnarray*}
\label{F2}
n \er^{0t}_n (\es^0_n)^{-1} \er^0_n & =&\bigg \{  \frac{1}{\theta^3} \sum  _ {l=1} ^ {k}  \ew_l^0  \big [  \big(  \frac{1}{n\theta}\sum  _ {i=1} ^ {k}  \ew_i^0 - \frac{1}{n (1-\theta )} \sum _{j =k+1} ^ n \ew_j^0 \big)^t \ew_l^0  \big]^2
\nonumber \\ && 
-\frac{1}{(1-\theta)^3}  \sum  _ {l =k+1} ^ n \ew_l^0  \big [ \big(\frac{1}{n\theta} \sum  _ {i=1} ^ {k}  \ew_i^0 - \frac{1}{n (1-\theta  )} \sum _{j =k+1} ^ n \ew_j^0 \big)^t \ew_l^0 \big]^2  
\nonumber \\ &&
\cdot\bigg (\frac{1}{n\theta^3} \sum  _ {l=1} ^ {k}  \ew_i^0  \big [\big( \frac{1}{n\theta}\sum  _ {i=1} ^ {k} \ew_i^0 - \frac{1}{n (1-\theta  )} \sum _{j =k+1} ^ n \ew_j^0 \big)^t\ew_l^0  \big]^2
\nonumber \\ && 
-\frac{1}{n(1-\theta)^3}  \sum  _ {l =k+1} ^ n \ew_l^0  \big [ \big( \frac{1}{n\theta} \sum  _ {i=1} ^ {k}   \ew_i^0 - \frac{1}{n (1-\theta  )} \sum _{j =k+1} ^ n \ew_j^0 \big)^t \ew_l^0 \big]^2 \bigg ) \bigg \} (1+o_{\eP}(1)).
\nonumber 
\end{eqnarray*}
Thus, using notations given by (\ref{eq8}) and (\ref{eq9}), we obtain
\begin{eqnarray*}
n \er^{0t}_n (\es^0_n)^{-1} \er^0_n &=& n \sum  _ {r,s,l,u,v=1} ^ {p} \omega ^r \omega ^s \big [  \frac{1}{n\theta^3} \sum  _ {i=1} ^ {k} w_{i,r}^0 w_{i,s}^0 w_{i,l}^0 -\frac{1}{n(1-\theta)^3} \sum  _ {j =k+1} ^ n  w_{j,r}^0 w_{j,s}^0 w_{j,l}^0 \big]     
\nonumber \\  
 && \cdot  \omega ^u \omega ^v  \big [ \frac{1}{n\theta^3} \sum  _ {i=1} ^ {k}  w_{i,u}^0 w_{i,v}^0 w_{i,l}^0 -\frac{1}{n(1-\theta)^3} \sum  _ {j =k+1} ^ n  w_{j,u}^0 w_{j,v}^0 w_{j,l}^0 \big] (1+o_{\eP}(1))\\
 & =&n \sum  _ {r,s,l,u,v=1} ^ {p} \alpha^{rsl} \alpha^{uvl} \omega ^r \omega ^s  \omega ^u \omega ^v(1+o_{\eP}(1)).
\end{eqnarray*}
In conclusion, for ${\cal E}_3$ of (\ref{F33}) and  for  $n \er^{0t}_n (\es^0_n)^{-1} \er^0_n$,  using assumptions (A6), (A7) and (A8), together with  the proof of Proposition 1 of  \cite{Guo:Zou:Wang:Chen:13}, we obtain:
$
{\cal E}_3=o_{\eP}(p^{1/2})$
and
$
n \er^{0t}_n (\es^0_n)^{-1} \er^0_n=o_{\eP}(p^{1/2})$.
 Combining the last two relations together relation (\ref{eq170}), we obtain that
\begin{equation*}
\label{eq141}
\EL_{nk}(\ebo)=n \epsi_n^0(\es^0_n)^{-1} \epsi_n ^0 +o_{\eP}(p^{1/2}).
\end{equation*}
\hspace*{\fill}$\blacksquare$ 
\ \\  

%%%%%%%%%%%%%%%%%%%%%%%%%%%%%%%%%%%%%%%%%%%%%%%%%%%%%%%%%%%%%%%%
%%%%%%%%%%%%%%%%%%%%%%%%%%%%%%%%%%%%%%%%%%%%%%%%%%%%%%%%%%%%%%%%
%%%%%%%%%%%%%%%%%%%%%%%%%%%%%%%%%%%%%%%%%%%%%%%%%%%%%%%%%%%%%%%%
 
\noindent {\bf Proof of Proposition \ref{proposition3}.} 
We first prove 
\begin{equation}
\label{eqp2}
n \epsi_{n}^{0t} \big((\ev^0_n)^{-1} -(\es^0_n)^{-1} \big )\epsi_n ^0= o_{\eP}(p^{1/2}).
\end{equation}
For this, we introduce the following two $p$-square matrices
\begin{equation*}
\label{Bn}
\ebb_n^0 \equiv (\ev_n^0)^{-1/2} \, \es_n^0 \, (\ev^0_n)^{-1/2}, \qquad    \ek_n^0\equiv \textbf{I}_p-\ebb_n^0
\end{equation*}
and the following  $p$-vector
\begin{equation*}
\label{mu}
 \eeta_n^0\equiv (\ev_n^0)^{-1/2} \epsi_n^0.
\end{equation*}
With this notations, the left hand side of relation (\ref{eqp2}), can be written
\begin{eqnarray*}
\label{eq33}
n \epsi_{n}^{0t} \big ((\ev_n^0)^{-1} -(\es_n^0)^{-1}\big )\epsi_n ^0
= n\eeta_n^{0t} \big ( \textbf{I}_p -(\ev_n^0)^{1/2} (\es^0_n)^{-1} (\ev_n^0)^{1/2}\big) \eeta_n^0
=  n\eeta_n^{0t} (\ek^0_n)^{-1} \eeta_n^0.
\end{eqnarray*}
We consider the following decomposition for  $ n\eeta_n^{0t} (\ek^0_n)^{-1} \eeta_n^0 $
\begin{equation}
\label{eq34}
 n\eeta_n^{0t} (\ek^0_n)^{-1} \eeta_n^0  
=\Big ( n\eeta_n^{0t} \ek_n^0 \eeta_n^0- n\eeta_n^{0t} (\ek^0_n)^2 \eeta_n^0-\cdots - (-1)^b n\eeta_n^{0t}(\ek^0_n)^b \eeta_n^0 \Big )
   +(-1)^b n \eeta_n^{0t} (\ek^0_n)^b ( \textbf{I}_p-(\ebb^0_n)^{-1} ) \eeta_n^0,
  \end{equation}
%\end{eqnarray}
for any $b \in \mathbb{N}^*$.\\
We will study the convergence of the expansion given by (\ref{eq34}). 
%For the term $n \eeta_n^{0t} \ek_n^0 \eeta_n^0$ of (\ref{eq34}), 
By Lemma 6 of \cite{Chen:Peng:Qin:09}, we have the inequality $\eeta_n^{0t} \textbf{A} \eeta_n^0  \leq \| \eeta_n^0\|^2 ( \tr(\textbf{A}^2))^{1/2}$, for any symmetric matrix $\textbf{A}$. Then, for the first term of the right-hand side of relation (\ref{eq34}) we have that, with probability one:
\begin{equation}
\label{eq38}
n\eeta_n^{0t} \ek_n^0 \eeta_n^0  \leq n \ \| \eeta_n^0\|^2 \big( \tr(\ek_n^0)^2\big)^{1/2}.
\end{equation}
Using assumption (A1) and Lemma 4 of \cite{Liu:Zou:Wang:13}, we obtain that, with probability one: 
\begin{equation*}
\label{eq35}
\|\eeta_n^0 \|^2 
= \epsi_n ^{0t} (\ev_n^0)^{-1}  \epsi_n ^0
\leq \frac{1}{\gamma_1(\ev_n^0)} \| \epsi_n^0\|^2
\leq \frac{1}{C_0}\| \, \epsi_n^0\|^2.
\end{equation*}
By Lemma \ref{lemma4}, we have that $\| \epsi_n^0\|= O_{\eP}(p^{1/2}n^{-1/2} )$ and thus
\begin{equation}
\label{eq36}
\|\eeta_n^0 \|^2=O_{\eP}(pn^{-1} ).
\end{equation}
By Lemma \ref{lemma7}, for $1 \leq r \leq p$, we have 
\begin{equation}
\label{eqqq}
\gamma_r(\ek^0_n) \leq (\tr (                                                                                                                                                                                                                                                                                                                                                                                                                                                                                                                                                                                                                                                                                                                                                                                                                                                                                                                                                                                                                                                                                                                                                                                                                                                                                                                                                                                                                                                                                                                                                                           \ek_n^0)^2)^{1/2}=O_{\eP}(p n^{-1/2}).
\end{equation}
Using relations (\ref{eq36}), (\ref{eqqq}) and  condition $p=o(n^{1/2})$ obtained by assumption (A6),  we have for (\ref{eq38}) that
\begin{equation}
\label{eqqq38}
n\eeta_n^{0t} \ek_n^0 \eeta_n^0  =   n \ O_{\eP}(n^{-1} p)  \ O_{\eP}(p n^{-1/2})= O_{\eP}(p^2  n^{-1/2})=o_{\eP}(p^{1/2}).
\end{equation}
On the other hand, using relations (\ref{eq36}) and (\ref{eqqq}), we obtain
\begin{equation*}
\label{eq37}
|\eeta_n^{0t} (\ek^0_n)^b \eeta_n^0| \leq  \|\eeta_n^0 \|^2 \max  _{1 \leq r \leq p} |\gamma_r(\ek^0_n)^b |
\leq \|\eeta_n^0 \|^2 (\tr(\ek^0_n)^2) ^{b/2}
=O_{\eP}(p\,n^{-1}) O_{\eP}(p^b n^{-b/2}),
\end{equation*}
which gives 
\begin{equation}
\label{eq43}
\eeta_n^{0t} (\ek^0_n)^b \eeta_n^0 =  O_{\eP}(p^{b+1} n^{-(b+2)/2}).
\end{equation}
The last equation means that the series $n \sum  _{b=1}^{\infty} (-1)^{b-1} \eeta_n^{0t} (\ek^0_n)^b \eeta_n^0 $ is convergent for fixed $n$ when $p=o(n^{1/2})$. Then, taking also into account relation  (\ref{eqqq38}), we can conclude that
\begin{equation}
\label{serie}
n \sum  _{b=1}^{\infty} (-1)^{b-1} \eeta_n^{0t} (\ek^0_n)^b \eeta_n^0 =o_{\eP}(p^{1/2}).
\end{equation}
The remaining task is to prove that the   term $n \eeta_n^{0t} (\ek_n^0)^b (\textbf{I}_p- (\ebb^0_n)^{-1})  \eeta_n^0$ in (\ref{eq34}) is negligible as $b \rightarrow \infty$.\\ 
For the last term of (\ref{eq34}), we have that
\begin{equation}
\label{eq41}
| n \eeta_n^{0t} (\ek_n^0)^b (\textbf{I}_p- (\ebb^0_n)^{-1})  \eeta_n^0|
\leq  | n \eeta_n^{0t} (\ek_n^0)^b   \eeta_n^0| +| n (\eeta_n^0)^t (\ek^0_n)^b  (\ebb^0_n)^{-1} \eeta_n^0|.
\end{equation} 
For the first term of the right hand side of (\ref{eq41}), by  relation (\ref{eq43}), we have, with probability one,  that 
$\eeta_n^{0t} (\ek^0_n)^b  \eeta_n^0 =  O_{\eP}(p^{b+1} n^{(-b-2)/2})$.
Then 
\begin{equation}
\label{eq45}
 n \eeta_n^{0t} (\ek^0_n)^b   \eeta_n^0=O_{\eP}(p^{b+1} n^{-b/2}).
\end{equation}
For the second term of the right hand side of (\ref{eq41}), we have, with probability one,  that
\begin{equation}
\label{eq46}
| n \eeta_n^{0t} (\ek^0_n)^b  (\ebb^0_n)^{-1} \eeta_n^0| \leq n \| \eeta_n^0 \|^2 \Mx((\ek^0_n)^b  (\ebb^0_n)^{-1}).
\end{equation}
Furthermore, according to Lemma 4 of \cite{Liu:Zou:Wang:13}, for any $p \times p$ symmetric matrix $\textbf{A}= (a_{ij})$, we have $\Mx(\textbf{A}) \leq \max_{1 \leq i \leq p} |\gamma_i(\textbf{A})|$. Then, with probability one, 
 \[\Mx((\ek^0_n)^b  (\ebb^0_n)^{-1}) \leq p \cdot \Mx (\ek^0_n)^b \cdot \Mx  (\ebb^0_n)^{-1} \leq p \cdot \max_{1 \leq r \leq p} |\gamma_r (\ek^0_n)^b|  \cdot  \max_{1 \leq r \leq p} |\gamma_r (\ebb^0_n)^{-1}|.\]
 By Lemma \ref{lemma7}, we know that $|\gamma_r(\ek^0_n)| \leq (\tr(\ek^0_n)^2)^{1/2}= O_{\eP}(p  n^{-1/2})$ and then,
 $|\gamma_r(\ek^0_n)^b| \leq (\tr(\ek^0_n)^2)^{b/2}= O_{\eP}(p^bn^{-b/2})$. On the other hand, by Proposition \ref{corolary2} it is clear that $\gamma_1(\ebb^0_n)^{-1}= \gamma_1((\es^0_n)^{-1} \ev^0_n) \leq \gamma_1(\ev^0_n) /C_0,$ with probability tending to one. All these imply  
  $\Mx((\ek^0_n)^b  (\ebb^0_n)^{-1}) \leq p   O_{\eP}(p^bn^{-b/2})$.
  Then, we obtain
\begin{equation}
\label{eq47}
\Mx((\ek^0_n)^b  (\ebb^0_n)^{-1})=O_{\eP}(p^{b+1}n^{-b/2}).
\end{equation}
Combining relations (\ref{eq36}), (\ref{eq46}) and (\ref{eq47}), for the second term of the right hand side of (\ref{eq41}), we obtain that
\begin{equation}
\label{eq48}
| n \eeta_n^{0t} (\ek^0_n)^b  (\ebb^0_n)^{-1} \eeta_n^0| \leq  n   p \| \eeta_n^0 \|^2 \Mx((\ek^0_n)^b  (\ebb^0_n)^{-1})
= O_p(p^{b+3}n^{-b/2}).
\end{equation}
By relations (\ref{eq41}), (\ref{eq45}), (\ref{eq48}) and assumption (A6), it follows that 
\begin{equation}
\label{eq49}
 n  \eeta_n^{0t} (\textbf{I}_p-\ev_n^0)  \eeta_n^0 = O_{\eP}(p^{b+1} n^{-b/2})+O_p(p^{b+3}n^{-b/2}) 
= O_{\eP}((p^{2+{4}/{b}} n^{-1})^b)
= o_{\eP}(1).
\end{equation}
Combining the results obtained in relations (\ref{serie}) and (\ref{eq49}),  we obtain (\ref{eqp2}). \\
The Proposition follows combining relation (\ref{eqp2}) and Proposition \ref{proposition2}. 
\hspace*{\fill}$\blacksquare$ 
\ \\ \ \\
%%%%%%%%%%%%%%%%%%%%%%%%%%%%%%%%%%%%%%%%%%%%%%%%%%%%%%%%%%%%%%%%
%%%%%%%%%%%%%%%%%%%%%%%%%%%%%%%%%%%%%%%%%%%%%%%%%%%%%%%%%%%%%%%%
%%%%%%%%%%%%%%%%%%%%%%%%%%%%%%%%%%%%%%%%%%%%%%%%%%%%%%%%%%%%%%%%

\noindent {\bf Proof of Theorem \ref{theorem2}.} 
 Since $\theta=k/n \rightarrow \theta^0 \in (0,1)$ and the point $k$, where the test is realised,  is known, we suppose that $\theta$ is $\theta^0$, then it is fixed.  \\
 
\textit{(i)} We prove  relation (\ref{eqt2}), by constructing a martingale and applying the martingale central limit theorem (see \cite{Chow:Teicher:97}). We will prove this relation in four steps. In Step 1 we construct a martingale, in Steps 2 and 3 we propose two sufficient conditions for applying a central limit theorem and finally, in Step 4 we prove  relation  (\ref{eqt2}). \\
%%%%%%%%%%%%%%%%%%%%%%%%%%%% Step 1
\hh \underline{\textit{Step 1.}} In this step, we will construct a martingale.\\
For $i =1, \cdots, n$, let us define the following  random vector sequence:
\[
\eg_i^0 \equiv  \left\{
\begin{array}{lll}
\displaystyle{ \frac{1}{\theta} \sum^i_{j=1} \ew^0_j}, & \textrm{if}& i \leq k, \\
 \displaystyle{ \frac{1}{\theta} \sum^k_{j=1} \ew^0_j - \frac{1}{1- \theta}\sum^i_{j=k+1} \ew^0_j}, & \textrm{if}& i > k,
\end{array}
\right.
\]
and also the random variable  $H_i^0 \equiv \|\eg_i^0\|^2-ip $. \\ 
Then, the left hand side of (\ref{eqt2}) can be written
 \begin{equation*}
 \label{eq56}
 \frac{n\epsi_{n}^{0t} (\ev_n^0)^{-1} \epsi_n^0  -p}{\Delta_n /n}= \frac{H^0_n}{\Delta_n}  .
 \end{equation*}
 The relation between $\eg_i^0$ and $\eg_{i-1}^0$ is:
 \[
 \eg_{i}^0=\left\{
 \begin{array}{lll}
\displaystyle{ \eg_{i-1}^0+\frac{ \ew_{i}^0}{\theta} }, & \textrm{ if } & i \leq k,\\
 & &\\
\displaystyle{ \eg_{i-1}^0- \frac{ \ew_i^0}{1-\theta}}, & \textrm{ if } & i >k,
 \end{array}
  \right.
 \]
with $\eg_0^0 \equiv \textbf{0}$.\\
Consider now the following filtration ${\cal F}_i=\sigma(\ew^0_1, \cdots , \ew^0_i)=\sigma(\eg_1^0, \cdots,\eg_i^0)$  for $ i=1\ldots,n$ the $\sigma$-field generated by $\ew^0_1, \cdots , \ew^0_i$ or by  $\eg_1^0, \cdots,\eg_i^0$. Firstly, we study if $\{H_i^0,{\cal F}_i\}_{i \geq 1}  $ is  a martingale.  For this, consider for example $i$ such that $i >k$. Then
\[
\eE[H^0_i | {\cal F}_{i-1} ] =  \|\eg^0_{i-1}\|^2 +\frac{\eE[\|\ew^0_i \|^2]}{(1- \theta)^2} - ip \neq H^0_{i-1}.
\]
  Consequently, $\{H_i^0,{\cal F}_{i}\}_{i \geq 1}  $ is not a martingale. We will now construct a martingale based on $\eg_n^0$ with respect to the filtration $\{{\cal F}_i\}_{i \geq 1}$.  For this, we define the following  random variable sequence $U_i \equiv H^0_i - H^0_{i-1}$, for $i =1, \cdots, n$, with $H_0^0=0$. Then
\[
 U_{i}^0=\left\{
 \begin{array}{lll}
\displaystyle{\frac{2}{\theta} \eg^{0t}_{i-1} \ew^0_i +  \frac{\|\ew^0_i \|^2}{\theta^2}   -p } , & \textrm{ if } & i \leq k,\\
 & &\\
\displaystyle{ -\frac{2}{1-\theta} \eg^{0t}_{i-1} \ew^0_i  +  \frac{\|\ew^0_i \|^2}{(1-\theta)^2}- p}, & \textrm{ if } &  i >k.
 \end{array}
  \right.
 \]  
  We consider the following two  random variable sequences:
\[\tau _i^0 \equiv U_i^0 -\eE[U_i^0],\]
 and 
  \[ 
 \varphi_i^0 \equiv \sum  _{j=1}^i \tau _j^0= H_i^0- \sum_{j=1}^i  \eE[U_j^0] .
 \]
 For all $i$ such that $i \leq k$ we have that   the condition expectation of $\varphi_i^0$  given the $\sigma$-field ${\cal F}_{i-1}$ is:
 \[
 \eE[ \varphi_i^0  | {\cal F}_{i-1} ]= \|\eg^0_{i-1}\|^2 +\frac{\eE[\|\ew^0_i \|^2]}{ \theta^2} - ip -\sum^{i}_{j=1} \eE[U^0_j]= \|\eg^0_{i-1}\|^2 -(i-1)p -\sum^{i-1}_{j=1} \eE[U^0_j]=\varphi_{i-1}^0
 \]
 and for all $i >k$:
  \[
 \eE[ \varphi_i^0  | {\cal F}_{i-1} ]= \|\eg^0_{i-1}\|^2 +\frac{\eE[\|\ew^0_i \|^2]}{ (1-\theta)^2} - ip -\sum^{i}_{j=1} \eE[U^0_j]= \varphi_{i-1}^0.
 \]
 Thus, $\{\varphi_i^0\}_{i \geq 1}  $ is  a martingale with respect to $\{{\cal F}_i\}_{i \geq 1}$.\\ 
 
To apply the martingale central limit theorem of \cite{Chow:Teicher:97} for $\{\varphi_i^0,{\cal F}_i, i \geq 1   \}$, it suffices to show that  
\begin{equation}
\label{ee}
\sum  _{i=1}^n \eE \big[ |\tau_i^0|^3\big] =o(\Delta^3_n)
\end{equation}
and 
\begin{equation}
\label{aa}
 \sum _{i=1}^n \eE \big[ |\eE\big [ (\tau_i^0)^2 | {\cal F}_{i-1} \big]- \sigma_i^2| \big]=o(\Delta^2_n),
\end{equation}
%%%%%%%%%%%%%%%%%%%%%%%%%%%%%%%%%% Step 2
\underline{\textit{Step 2.}} In this step, we will prove relation (\ref{ee}). \\
In order to facilitate writing, for $i=1, \cdots,n$, we denote
\begin{equation}
\label{N}
\en_i^0 \equiv \frac{\ew_{i}^0}{\theta} \e1_{i \leq k} -\frac{\ew_i^0}{1-\theta} \e1_{i>k}.
\end{equation} 
Then, the  random variable $\tau _i^0$ can be written  
\begin{eqnarray}
\label{eq67}
\tau _i^0 = 2 \eg_{i-1}^{0t}  \en_i^0 + \|\en_i^0 \|^2 -\eE [\|\en_i ^0\|^2].
\end{eqnarray}
By assumption (A9), for some positive absolute constant $C_7<\infty$, for all $ i=1,\cdots, n$ and all $j_1,\ldots, j_l=1, \ldots, p$, $l \in \mathbb{N}$, whenever $\sum  _{i=1} ^l d_i \leq 6$, we have 
\begin{equation}
\label{C7}
\eE[ N_{i,j_1} ^{ 0d_1} \cdots  N_{i,j_l} ^{0d_l}] \leq C_7,
 \end{equation}
 with $N_{i,j_1} ^0$ is the $j_1$-th components of the vector $\en_i^0$ defined in (\ref{N}).
By the Holder inequality, for any $b \leq 3$ we have
\begin{equation}
\label{eqni}
\eE [ \|\en_i^0\|^{2b} ]\leq C \eE[ \big(\sum^p_{j=1} w^2_{ij}\big)^b] \leq C p^{b-1} \sum^p_{j=1} \eE[w^{2b}_{ij} ]\leq C p^b.
\end{equation}
By the Cauchy Schwartz's inequality, we have that
\begin{equation}
\label{eq68}
\eE\big[ \en_i ^{0t} \eg_{i-1}^0   \big]^3 \leq \big[\eE [  \en_i ^{0t} \eg_{i-1} ^0  ]^6\big ]^{1/2}.
\end{equation}
On the other hand, by Lemma 7 of \cite{Guo:Zou:Wang:Chen:13}, we have
\begin{equation}
\label{eq69}
\eE[ \en_i ^{0t} \eg_{i-1} ^0  ]^6 \leq C\,p^6(i^3+i^2+i).
\end{equation}
Then, by (\ref{eq68}) and (\ref{eq69}), we obtain that  
$\eE \big[ \en_i ^{0t}\eg_{i-1}^0   \big]^3  \leq [C\,p^6(i^3+i^2+i)]^{1/2}$. 
This implies that 
\begin{equation}
\label{eq70}
\eE \big[ \en_i ^{0t} \eg_{i-1} ^0 \big ]^3 \leq C\,p^3(i^{3/2}+i+i^{1/2}).
\end{equation}
On the other hand, we have
\begin{equation*}
\label{eq71} \eE \big [\|\en_i^0 \|^2 -\eE [\|\en_i^0 \|^2 ]\big]^3 
=\eE \left[\|\en_i^0 \|^6 - 3\|\en_i^0 \|^4 \eE [\|\en_i ^0\|^2] +3 \|\en_i^0 \|^2 (\eE [\|\en_i^0 \|^2])^2-(\eE [\|\en_i^0 \|^2])^3\right ].
\end{equation*}
Using (\ref{eqni}), we obtain that 
\begin{equation}
\label{eq72}
\left| \eE \big[\|\en_i^0 \|^2 -\eE [\|\en_i ^0\|^2] \big]^3 \right| \leq Cp^3.
\end{equation}
Using relations (\ref{eq70}) and (\ref{eq72}), we can write
$
\eE [|\tau_i^0|^3] \leq Cp^3 (1+i^{3/2}+i+i^{1/2}).
$ 
Then 
$
\sum_{i=1}^n \eE [|\tau_i^0|^3] \leq  Cnp^3 (1+n^{3/2}+n+n^{1/2}),
$ 
which gives 
\begin{equation}
\label{eq75}
\sum  _{i=1}^n \eE [|\tau_i^0|^3] \leq Cp^3 (n^{5/2}+n^2+n^{3/2} +n).
\end{equation}
On the other hand, using assumption (A1), by similar arguments as for  relation (16) of \cite{Guo:Zou:Wang:Chen:13}, we have that $ \Delta^2_n \geq C n^2 p +O_p(p^2n)$. Then
\begin{equation}
\label{eq76}
 \Delta^3_n \geq Cn^3p^{3/2} +O_p(p^3n^{3/2}).
\end{equation}
From relations (\ref{eq75}) and (\ref{eq76}), we obtain 
\begin{eqnarray*}
\label{eq77}
\frac{\sum _{i=1}^n \eE [|\tau_i^0|^3]}{\Delta^3_n}  \leq C \,  \frac{p^3 n^{5/2}}{n^3 p^{3/2}} = O_p \big((p^3/n) ^{1/2} \big).
\end{eqnarray*}
Since $p=o(n^{1/3})$,  relation (\ref{ee}) follows. \\ 
%%%%%%%%%%%%%%%%%%%%%%%%%%%%%%%%%% Step 3
\hh \underline{\textit{Step 3.}} Now, in this step we prove relation (\ref{aa}). \\
By  elementary calculations and  using relation (\ref{eqni}), we obtain
$(\tau_i^0)^2= 4 \eg_{i-1}^{0t}  \en_i^0  \en_i ^{0t} \eg_{i-1}^0 + 4 \eg_{i-1}^{0t}  \en_i^0 (\|\en_i^0\|^2-\eE [\|\en_i^0\|^2])+O_{\eP}(p^2)$. We observe also that $\sigma_i^2=\eE \big[(\tau_i^0)^2 \big]$.
Then 
\begin{equation}
\label{eq80}
\eE \big[(\tau_i^0)^2 | {\cal F}_{i-1} \big] = 4 \eg_{i-1}^{0t}\eE \big[ \en_i^0  \en_i ^{0t} \big]\eg_{i-1}^0+ 4 \eg_{i-1}^{0t}  \eE\big [ \en_i^0 \big(\|\en_i^0\|^2-\eE [\|\en_i^0\|^2] \big )\big] +O_{\eP}(p^2)
\end{equation}
and 
\[
\sigma_i^2=4 \eE \big [ \eg_{i-1}^{0t}\eE[ \en_i^0  \en_i ^{0t}] \eg_{i-1}^0  \big ]+\eE\big[ \big(\|\en_i^0\|^2-\eE [\|\en_i^0\|^2 ] \big)^2  \big] +O(p^2).
\]
By inequality (\ref{eqni}) we have $\eE\big[ \big(\|\en_i^0\|^2-\eE [\|\en_i^0\|^2 ] \big)^2  \big]=O(p^2)$. Then 
\begin{equation}
\label{sig2}
\sigma_i^2=4 \eE \big [ \eg_{i-1}^{0t}\eE[ \en_i^0  \en_i ^{0t}] \eg_{i-1}^0  \big ] +O(p^2).
\end{equation}
Using relations (\ref{eq80}) and (\ref{sig2}), we obtain that, for any $i=1, \cdots, n$, 
\begin{eqnarray}
\label{eq82}
  \eE \big [ \big|\eE [\tau_i^2 |{\cal F}_{i-1} ]- \sigma_i^2 \big|^2 \big ]
&\leq &  16 \bigg ( \eE  \big [ \eg_{i-1}^{0t} \eE[ \en_i^0  \en_i ^{0t}] \eg_{i-1}^0- \eE [\eg_{i-1}^{0t} \en_i^0  \en_i ^{0t} \eg_{i-1}^0 ]  \big]  ^2
\nonumber \\ &&
+ \eE \big [ \en_i^{0t}  \big (\|\en_i^0\|^2-\eE [\|\en_i^0\|^2 ]\big)  \big ]
\eE \big[\eg_{i-1}^0 \eg_{i-1}^{0t} \big] 
\eE \big[\en_i^0  \big(\|\en_i^0\|^2-\eE [\|\en_i^0\|^2] \big) \big]
+O(p^4)\bigg ) 
\nonumber \\  
&\equiv& 16 \big ( A+B+O(p^4) \big ) 
\end{eqnarray}
\underline{For the term A of (\ref{eq82})}, we have the decomposition
\begin{equation}
\label{eq83}
A= \eE \big[\eg_{i-1}^{0t} \eE[ \en_i^0  \en_i ^{0t}]\eg_{i-1}^0\big]^2 
-  \big(\eE[\eg_{i-1}^{0t}  \en_i^0 \en_i ^{0t} \eg_{i-1}^0] \big)^2
 \equiv  A_1- A_2.
\end{equation}
Before analysing the terms $A_1$ and $A_2$, we note that
\begin{equation}
\label{eq84}
\eE(\en_i^0  \en_i ^{0t})
=(\ev_n^0)^{-1/2}  \big[\frac{1}{\theta^2} \e1_{i \leq k}+\frac{1}{(1-\theta)^2} \e1_{i>k}\big] \ev_{(i)}^0 (\ev_n^0)^{-1/2} .
\end{equation}
In order to facilitate writing, we consider the following matrix
\[\em_i^0 \equiv \left( \frac{1}{\theta^2} \e1_{i \leq k}+\frac{1}{(1-\theta)^2} \e1_{i>k}  \right) \ev_{(i)}^0.\] Then,   $\eE[\en_i^0  \en_i ^{0t}]$ can be expressed as:
\begin{equation}
\label{eq85}
\eE[\en_i^0  \en_i ^{0t}]=(\ev_n^0)^{-1/2} \, \em_i^0 \, (\ev_n^0)^{-1/2}.
\end{equation}
Hence, the term $A_1$ of (\ref{eq83}), can be written
\begin{equation*}
A_1 =
\eE \big [\eg_{i-1}^{0t}       (\ev_n^0)^{-1/2}  \em_i^0 \ (\ev_n^0)^{-1/2} 
\eg_{i-1}^{0}   \big] ^2.
\end{equation*}
Taking into account  that the random vectors $\ew^0_i$ are independent, with mean zero, for all $ i=1,\ldots, n$, we can decompose $A_1$ as 
\begin{equation}
\label{eq86}
A_1\equiv  A_{11} + A_{12} + A_{13},
\end{equation} 
with \\
\textit{ - if } $i-1 \leq k$:
\[
A_{11}  \equiv  \frac{1}{\theta^4} \sum^{i-1}_{h \neq l =1} \eE \big[\ew_h ^0 (\ev_n^0)^{-1/2}  \em_i^0 \ (\ev_n^0)^{-1/2}  \ew_h ^0 \big] \eE \big[\ew_l ^0 (\ev_n^0)^{-1/2}  \em_i^0 \ (\ev_n^0)^{-1/2}  \ew_l ^0 \big],
\]
\[
 A_{12} \equiv \frac{2}{\theta^4}  \sum  _{h\neq l=1}^{i-1} \eE \big[\ew^0_h (\ev_n^0)^{-1/2}  \em_i^0 \, (\ev_n^0)^{-1/2} \ew^0_l \ew^{0t}_h (\ev_n^0)^{-1/2}  \em_i^0 \, (\ev_n^0)^{-1/2} \ew^0_l \big]
\]
\[
A_{13} \equiv   \frac{1}{\theta^4}   \sum  _{h=1}^{i-1}\eE \big  [ \ew_h ^{0t} \  (\ev_n^0)^{-1/2}  \em_i^0 \ (\ev_n^0)^{-1/2}  \  \ew_h ^0 \big  ]^2 .
\]
\textit{ - if } $i-1 > k$:
\begin{eqnarray*}
A_{11} & \equiv &  \frac{1}{\theta^4} \sum^{k}_{h \neq l =1} \eE \big[\ew_h ^0 (\ev_n^0)^{-1/2}  \em_i^0 \ (\ev_n^0)^{-1/2}  \ew_h ^0 \big] \eE \big[\ew_l ^0 (\ev_n^0)^{-1/2}  \em_i^0 \ (\ev_n^0)^{-1/2}  \ew_l ^0 \big]
\nonumber \\
 && + \frac{1}{(1-\theta)^4}  \sum^{i-1 }_{h' \neq l' =k+1} \eE \big[\ew_{h'} ^0 (\ev_n^0)^{-1/2}  \em_i^0 \ (\ev_n^0)^{-1/2}  \ew_{h'} ^0 \big] \eE \big[\ew_{l'} ^0 (\ev_n^0)^{-1/2}  \em_i^0 \ (\ev_n^0)^{-1/2}  \ew_{l'} ^0 \big],
\nonumber 
\end{eqnarray*}
%%%%%%%%%%%%%%%%%%%%%%%%%%%%%%%%%%%%
\begin{eqnarray*}
 A_{12}& \equiv & \frac{2}{\theta^4}  \sum  _{h\neq l=1}^{k} \eE \big[\ew^0_h (\ev_n^0)^{-1/2}  \em_i^0 \, (\ev_n^0)^{-1/2} \ew^0_l \ew^{0t}_h (\ev_n^0)^{-1/2}  \em_i^0 \, (\ev_n^0)^{-1/2} \ew^0_l \big]
 \nonumber \\  
 &+& \frac{2}{(1-\theta)^4}   \sum  _{h'\neq l'=k+1}^{i-1}\eE \big[ \ew^0_{h'} (\ev_n^0)^{-1/2}  \em_i^0 \, (\ev_n^0)^{-1/2} \ew^0_{l'} \ew^{0t}_{h'} (\ev_n^0)^{-1/2}  \em_i^0 \, (\ev_n^0)^{-1/2} \ew^0_{l'}\big]. \\
  \nonumber 
  \end{eqnarray*}
%%%%%%%%%%%%%%%%%%%%%%%%%%%%%%%%%%%%
\begin{eqnarray*}
A_{13}& \equiv &  \frac{1}{\theta^4}   \sum  _{h=1}^{k}\eE \big  [ \ew_h ^{0t} \  (\ev_n^0)^{-1/2}  \em_i^0 \ (\ev_n^0)^{-1/2}  \  \ew_h ^0 \big  ]^2 
%\nonumber \\&& 
+  \frac{1}{(1-\theta)^4} \sum  _{h'=k+1}^{i-1} \eE \big  [  \ew_{h'} ^{0t}  \ (\ev_n^0)^{-1/2}  \em_i^0 \ (\ev_n^0)^{-1/2} \  \ew_{h'} ^0 \big  ]^2.
\nonumber 
\end{eqnarray*}
We study $ A_{11}$, $ A_{12}$, $ A_{13}$. For this, we consider the case $i-1 > k$, the other is similar.  \\
For $ A_{11}$, applying  Cauchy-Schwarz  inequality, we have
\begin{eqnarray*}
\label{eq87}
A_{11}& \leq &  \frac{1}{\theta^4}  \Big  (  \sum  _{h=1}^{k} \eE \big [ \ew_h ^{0t}   (\ev_n^0)^{-1/2}  \em_i^0 \ (\ev_n^0)^{-1/2}  \ew_h ^0 \big ] \Big  )^2 
 + \frac{1}{(1-\theta)^4} \Big  (   \sum  _{h'=k+1}^{i-1} \eE \big [ \ew_{h'} ^{0t}   (\ev_n^0)^{-1/2}  \em_i^0 \ (\ev_n^0)^{-1/2} \ew_{h'} ^0  \big ]\Big  )^2
\nonumber \\ 
&=&  \frac{1}{\theta^4} \Big  ( \sum  _{h=1}^{k}  \eE \big [  \ez_h ^{0t}   (\ev_n^0)^{-1}  \em_i^0 \ (\ev_n^0)^{-1}  \ez_h ^0 \big ] \Big  )^2 
 +  \frac{1}{(1-\theta)^4} \Big  (  \sum  _{h'=k+1}^{i-1} \eE \big [ \ez_{h'} ^{0t}   (\ev_n^0)^{-1}  \em_i^0 \ (\ev_n^0)^{-1} \ez_{h'} ^0  \big ]\Big  )^2.
\end{eqnarray*}
Then,  since $(\varepsilon_i)_{1 \leq i \leq n}$ are independent, we have:
\begin{equation}
\label{eq88}
A_{11}= \frac{1}{\theta^4}\Big  ( \sum  _{h=1}^{k} \tr\big(\ev_{(h)}^0  (\ev_n^0)^{-1}  \em_i^0 (\ev_n^0)^{-1}\big) \Big )^2 + \frac{1}{(1-\theta)^4}  \Big  ( \sum _{h'=k+1}^{i-1} \tr \big (\ev_{(h')}^0  (\ev_n^0)^{-1}  \em_i^0 \  (\ev_n^0)^{-1} \big)
\Big )^2.
\end{equation}
Similarly, we have for $ A_{12}$:
\begin{eqnarray*}
\label{eq89}
 A_{12}&=& \frac{2}{\theta^4}  \eE \Big[ \sum  _{h\neq l=1}^{k} \ew^0_h (\ev_n^0)^{-1/2}  \em_i^0  (\ev_n^0)^{-1/2} \ew^0_l \ew^{0t}_l (\ev_n^0)^{-1/2}  \em_i^0  (\ev_n^0)^{-1/2} \ew^0_h \Big]
 \nonumber \\ && 
 + \frac{2}{(1-\theta)^4}  \eE \Big[ \sum  _{h'\neq l'=k+1}^{i-1} \ew^0_{h'} (\ev_n^0)^{-1/2}  \em_i^0  (\ev_n^0)^{-1/2} \ew^0_{l'} \ew^{0t}_{l'} (\ev_n^0)^{-1/2}  \em_i^0  (\ev_n^0)^{-1/2} \ew^0_{h'}\Big]
  \nonumber \\ 
 &=& \frac{2}{\theta^4}  \eE \Big[ \sum  _{h\neq l=1}^{k} \ew^0_h   (\ev_n^0)^{-1/2}  \em_i^0  (\ev_n^0)^{-1/2}  \eE \big[  \ew^0_l \ew^{0t}_l \big]   (\ev_n^0)^{-1/2}  \em_i^0  (\ev_n^0)^{-1/2} \ew^0_h \Big]
 \nonumber \\ 
 && 
 + \frac{2}{(1-\theta)^4}  \eE \Big[ \sum  _{h'\neq l'=k+1}^{i-1} \ew^0_{h'} (\ev_n^0)^{-1/2}  \em_i^0  (\ev_n^0)^{-1/2} \eE \big[ \ew^0_{l'} \ew^{0t}_{l'} \big] (\ev_n^0)^{-1/2}  \em_i^0  (\ev_n^0)^{-1/2} \ew^0_{h'}\Big]
  \nonumber \\
  &\leq & \frac{2}{\theta^4} p  k^2   \gamma_1^4(\ev_n^0)^{-1}  
 \sup_{h,l \in \{1, \cdots, k\}} \big(\gamma_1(\ev_{(l)}^0 ) \, \gamma_1 (\ev_{(h)}^0)\big) \cdot \gamma_1^2 (\em_i ^0)
   \nonumber \\ 
 && 
 + \frac{2}{(1-\theta)^4}  p \big(i-1-k\big)^2 \, \gamma_1^4(\ev_n^0)^{-1}  \sup_{h',l' \in \{k+1, \cdots, i-1\}} \big(\gamma_1(\ev_{(l')}^0) \, \gamma_1(\ev_{(h')}^0)\big) \cdot \gamma_1^2 (\em_i ^0).
 \nonumber \\
  \end{eqnarray*}
Taking into account  assumption (A1), we obtain that
\begin{equation}
\label{eq90}
A_{12} \leq Cp (i-1)^2.
\end{equation}
For the term $A_{13}$ of (\ref{eq86}), we have
\begin{eqnarray}
A_{13} &=& \frac{1}{\theta^4} \sum  _{h=1}^{k} \, \sum  _{j,l,s,t=1}^{p} \Big( (\ev_n^0)^{-1/2}  \em_i^0  (\ev_n^0)^{-1/2} \Big)_{jl}  \Big( (\ev_n^0)^{-1/2}   \em_i^0  (\ev_n^0)^{-1/2} \Big)_{st} \eE[w_{h,j}^0 w_{h,l} ^0 w_{h,s}^0  w_{h,t}^0]
 \nonumber \\ && 
 + \frac{1}{(1-\theta)^4} \sum  _{h'=k+1}^{i-1} \, \sum  _{j,l,s,t=1}^{p} \Big( (\ev_n^0)^{-1/2}  \em_i^0  (\ev_n^0)^{-1/2} \Big)_{jl}  \Big( (\ev_n^0)^{-1/2}   \em_i^0   (\ev_n^0)^{-1/2}\Big )_{st} \eE[w_{h',j}^0 w_{h',l}^0 w_{h',s} ^0 w_{h',t}^0]. \nonumber
\end{eqnarray}
Taking into account  assumption (A9), we obtain that
\begin{eqnarray*}
\label{eq94}
 A_{13}&\leq& \frac{1}{\theta^4} C  p^4 k  \gamma_1^4 (\ev_n^0)^{-1}    \gamma_1^2 (\em_i ^0) +  \frac{1}{(1-\theta)^4} C  p^4 (i-k)  \gamma_1^4 (\ev_n^0)^{-1}    \gamma_1^2 (\em _i^0).
\end{eqnarray*}
Using also  assumption (A1), we obtain get
\begin{equation}
\label{eq95}
A_{13} \leq Cp^4 (i-1).
\end{equation}
For the term $A_2$  of (\ref{eq83}), by  similar calculations,  we obtain  that\\
\textit{ - if } $i-1 \leq k$:
\[
A_2=\frac{1}{\theta^4}  \Big  (\sum  _{h=1}^{i-1} \tr \big(\ev_{(h)}^0  (\ev_n^0)^{-1}  \em_i^0 (\ev_n^0)^{-1}\big) \Big )^2 
\]
\textit{ - if } $i-1 > k$:
\begin{equation}
\label{eq96}
A_{2} = \frac{1}{\theta^4}  \Big  (\sum  _{h=1}^{k} \tr\big(\ev_{(h)}^0  (\ev_n^0)^{-1}  \em_i^0 (\ev_n^0)^{-1} \big) \Big )^2 + \frac{1}{(1-\theta)^4} \Big  (  \sum _{h'=k+1}^{i-1} \tr \big (\ev_{(h')}^0  (\ev_n^0)^{-1}  \em_i^0 \  (\ev_n^0)^{-1} \big)
\Big )^2.
\end{equation}
For $i-1 > k$, by relations (\ref{eq88}) and  (\ref{eq96}), we obtain that 
\begin{equation}
\label{eq99}
 A_{11}-A_2 \leq 0.
\end{equation}
Inequality (\ref{eq99})  is also true for $i-1 \leq k$. \\
Then, since the term A of (\ref{eq83}) can be written $A=(A_{11}+A_{12}+A_{13})- A_2$, combining relations (\ref{eq90}), (\ref{eq95}) and (\ref{eq99}), we obtain that  
\begin{equation}
\label{eq100}
 A \leq Cp(i-1)^2+ Cp^4 (i-1).
\end{equation}
\underline{For the term B of (\ref{eq82})}, in the case  $i-1 > k$,  taking into account the fact that $\eE(\en_h^0)= \textbf{0}_p$, $\eE(\ew_h^0)= \textbf{0}_p$ for all $h=1, \ldots,n$ and the fact $\ew_h^0$ is independent of $\ew_{h'}^0$ for $h \neq h'$, we have that
\begin{eqnarray*}
\label{eq101}
B&=& \eE \big [\en_i^{0t} \| \en_i^0\|^2 \big]  \eE \big [ \big (\frac{1}{\theta} \sum  _{h=1}^{k} \ew_h^0 - \frac{1}{1-\theta} \sum _{h'=k+1}^{i-1}  \ew_{h'} ^0 \big)
%\nonumber \\ &&
\cdot \big(\frac{1}{\theta} \sum  _{h=1}^{k} \ew_h - \frac{1}{1-\theta} \sum _{h'=k+1}^{i-1}  \ew_{h'} \big)^t \big ] \eE \big [\en_i^{0} \| \en_i^0\|^2 \big]
\nonumber \\
&=& \frac{1}{\theta^2} \eE \big[\en_i^{0t} \| \en_i^0\|^2 \big]  \sum  _{h,l=1}^{k} \eE \big[\ew_h^0  \ew_l^{0t} \big] \, \eE \big[\en_i^{0t} \| \en_i^0\|^2 \big]
%\nonumber \\ &&
 + \frac{1}{(1-\theta)^2}\eE \big  [\en_i^{0t} \| \en_i^0\|^2 \big ]  \sum  _{h',l'=k+1}^{i-1}  \eE \big [\ew_{h'}^0  \ew_{l'}^{0t} \big ] \, \eE \big [ \en_i^{0} \| \en_i^0\|^2 \big ] 
\nonumber \\ 
&=& \frac{1}{\theta^2}\sum  _{u,s=1}^{p}  \sum  _{l,r=1}^{p} \Big (\sum  _{h=1}^{k} \ev_n^{-1/2} \ev_{(h)}^0 \ev_n^{-1/2} \Big )_{us}\, \eE \big[N_{i,u}^0 (N_{i,l}^0)^2\big ] \eE \big[N_{i,s} ^0 (N_{i,r}^0)^2 \big]
\nonumber \\ &&
+\frac{1}{(1-\theta)^2} \sum  _{u,s=1}^{p} \sum  _{l,r=1}^{p} \Big ( \sum  _{h'=k+1}^{i-1}  (\ev_n^0)^{-1/2}  \ev_{(h')}^0 (\ev_n^0)^{-\frac{1}{2}}  \Big )_{us} 
\, \eE \big[N_{i,u}^0 (N_{i,l}^0)^2\big ] \eE \big[N_{i,s} ^0 (N_{i,r}^0)^2 \big ],
\end{eqnarray*}
where $\Big (\sum  _{h=1}^{k} \ev_n^{-1/2} \ev_{(h)}^0 \ev_n^{-1/2} \Big )_{us}$ is the $(u,s)$-th element of the matrix $\sum  _{h=1}^{k} \ev_n^{-1/2} \ev_{(h)}^0 \ev_n^{-1/2}$.
Using Lemma 4 of \cite{Liu:Zou:Wang:13} and relation (\ref{C7}), we get
\begin{eqnarray*}
B&\leq& \frac{1}{\theta^2} \sum  _{u,s=1}^{p}  \sum  _{l,r=1}^{p} \Mx \Big ( \sum  _{h=1}^{k} \ev_n^{-1/2} \ev_{(h)}^0 \ev_n^{-1/2}  \Big )_{us}  \eE \big[N_{i,u}^0 (N_{i,l}^0)^2\big ] \eE \big[N_{i,s} ^0 (N_{i,r}^0)^2 \big]
\nonumber \\ &&
+\frac{1}{(1-\theta)^2} \sum  _{u,s=1}^{p} \sum  _{l,r=1}^{p} \Mx \Big ( \sum  _{h'=k+1}^{i-1}  (\ev_n^0)^{-1/2}  \ev_{(h')}^0 (\ev_n^0)^{-1/2}  \Big )_{us}   \eE \big[N_{i,u}^0 (N_{i,l}^0)^2\big ] \eE \big[N_{i,s} ^0 (N_{i,r}^0)^2 \big ]
\nonumber \\ 
& \leq  &\frac{ C}{\theta^2} p^4 (k \gamma_1^2(\ev_n^{-1/2}) \sup_{h \in \{1, \cdots, k \}} \big( \gamma_1(\ev_{(h)}^0 ) \big) +\frac{ C}{(1-\theta)^2} p^4 (i-1-k)\gamma_1^2(\ev_n^{-1/2})   \sup_{h \in \{k+1, \cdots, i-1 \}} \big(\gamma_1(\ev_{(h')} ^0 ) \big).
\end{eqnarray*}
Then, by assumption (A1), we obtain that
\begin{equation}
\label{eq102}
 B \leq Cp^4(i-1).
\end{equation} 
Similarly, we can prove that  inequality (\ref{eq102})  is also true for $i-1 \leq k$. \\
In conclusion, combining relations (\ref{eq82}), (\ref{eq100}) and (\ref{eq102}), we get that
\begin{equation*}
\label{eq103}
 \eE \big [ |\eE[(\tau_i^0)^2 | {\cal F}_{i-1}]- \sigma_i^2|^2 \big ]
\leq   Cp(i-1)^2 + Cp^4(i-1) +O(p^4) .
\end{equation*}
Then, by Cauchy-Schwarz inequality, we have
\begin{eqnarray*}
\label{eq104}
\frac{\sum _{i=1}^n \eE \big [ |\eE [(\tau_i^0)^2 | {\cal F}_{i-1}]- \sigma_i^2| \big ] }{\Delta^2_n} &\leq & \frac{\sum _{i=1}^n \Big (\eE\big[|\eE[(\tau_i^0)^2 | {\cal F}_{i-1}]- \sigma_i^2|^2 \big] \Big)^{1/2} }{\Delta^2_n}
\nonumber \\ 
 &\leq& \frac{C n (n^{1/2} p^2 + p^{1/2} n + O(p^2))}{C \, n^2 p+ O_{\eP}(p^2 n)}
\nonumber \\ 
 &\leq& C \frac{n^{3/2} p^2}{n^2 p}.
\end{eqnarray*}
and hence the relation (\ref{aa}) follows. \\
%%%%%%%%%%%%%%%%%%%%%%%%% Step 4
\hh \underline{\textit{Step 4.}} In this step, on the basis of the central limit theorem for martingales, we will complete the proof of  relation  (\ref{eqt2}). On the basis of relations (\ref{ee}) and (\ref{aa}) proved in Step 2 and Step 3, applying the martingale central limit Theorem of \cite{Chow:Teicher:97} (Theorem 1, page 336), for $\{\varphi_n^0, {\cal F}_n\}_{n \geq 1}$ we get:
\begin{equation}
\label{fin1}
 \frac{\varphi_n^0}{\Delta_n}\overset{{\cal L}} {\underset{n \rightarrow \infty}{\longrightarrow}}  { \cal N}(0,1).
 \end{equation}
Using notations given in \textit{Step 1}, we have 
$n \left( n\epsi_{n} ^{0t} (\ev_n^0)^{-1} \epsi_n^0  -p\right)=  H^0_n = \varphi_n^0 +\sum _{i=1}^n  \eE[U_i^0].$
In the other hand, since $H_0^0=0$, we have $\sum  _{i=1}^n  \eE[U_i^0] =\eE[H^0_n]=\eE[\| \eg_n^0 \|^2]-np$. But, taking into account relations (\ref{eq4}) and (\ref{var}), we get:
\begin{eqnarray*}
\eE[\| \eg_n^0 \|^2] & = & \frac{1}{\theta^2} \sum^k_{i=1} \eE[\| \ew_i^0 \|^2] + \frac{1}{(1-\theta)^2} \sum^n_{j=k+1} \eE[\| \ew_j^0 \|^2]   \\
& = & \frac{\sigma^2}{\theta^2} \sum^k_{i=1} \big( \XX_i^t (\ev_n^0)^{-1/2} \big) \big((\ev_n^0)^{-1/2} \XX_i \big) +\frac{\sigma^2}{(1-\theta)^2} \sum^n_{j=k+1} \big( \XX_j^t (\ev_n^0)^{-1/2} \big) \big((\ev_n^0)^{-1/2} \XX_j \big) \\
& = & \frac{\sigma^2}{\theta^2} \sum^k_{i=1} \tr \bigg( \big( \XX_i^t (\ev_n^0)^{-1/2} \big) \big((\ev_n^0)^{-1/2} \XX_i \big) \bigg) +\frac{\sigma^2}{(1-\theta)^2} \sum^n_{j=k+1}  \tr \bigg( \big( \XX_j^t (\ev_n^0)^{-1/2} \big) \big((\ev_n^0)^{-1/2} \XX_j \big) \bigg) \\
&=& \sigma^2 \tr \big( (\ev_n^0)^{-1/2}\big( \frac{1}{\theta^2} \sum^k_{i=1} \XX_i \XX_i^t + \frac{1}{(1-\theta)^2} \sum^n_{j=k+1} \XX_j \XX_j^t \big)  (\ev_n^0)^{-1/2}  \big) \\
& = & n \tr \big( (\ev_n^0)^{-1/2} \ev_n^0 (\ev_n^0)^{-1/2} \big) \\
&= & np.
\end{eqnarray*}
Thus $\eE[H_n^0]=0$. 
Then, taking into account relation (\ref{fin1}), we obtain  claim (\ref{eqt2}).\\

\textit{(ii)} The assertion results from (i) combined with  Proposition \ref{proposition2} and Proposition \ref{proposition3}.
\hspace*{\fill}$\blacksquare$ \\

%%%%%%%%%%%%%%%%%%%%%%%%%%%%%%%%%%%%%%%%%ù
\noindent {\bf Proof of Theorem  \ref{theorem_H1}.} 
By elementary calculations, we have under hypothesis $H_1$, with probability one:
\[
\epsi_n (\ebo)=\epsi^0_n - \frac{1}{n(1-\theta)} \bigg(\sum^n_{j=k+1}\XX_j \XX_j^t  \bigg)(\eb^0_2 -\ebo).
\]
Also, under $H_1$, since $\eE[\ez_i^0]=\textbf{0}_p$, using   assumption (A1), we have
\[
\epsi_n (\ebo)= - \frac{1}{n(1-\theta)} \bigg(\sum^n_{j=k+1} \XX_j \XX_j^t \bigg) (\eb^0_2 -\ebo)(1+o_{\eP}(1)).
\]
Then, 
\[
n\epsi _n^{t}(\ebo) \big(\ev_n^0\big)^{-1} \epsi_n(\ebo)=\Big(n \epsi _n^{0t} \big(\ev_n^0\big)^{-1} \epsi _n^{0} +n \big| (\eb^0_2 -\ebo)^t \frac{1}{n(1-\theta)} \sum^n_{j=k+1} \XX_j \XX_j^t  \big|^2 \Big)(1+o_{\eP}(1)).
\]
The test statistic becomes
\[
{\cal Z}(\ebo)=\frac{n \epsi _n^{0t} \big(\ev_n^0\big)^{-1} \epsi _n^{0} -p}{\Delta_n /n}+ n \frac{\big| (\eb^0_2 -\ebo)^t \frac{1}{n(1-\theta)} \sum^n_{j=k+1} \XX_j \XX_j^t  \big|^2}{\Delta_n /n} (1+o_{\eP}(1)).
\]
Since $\eb^0_2 \neq \ebo$, using with Theorem \ref{theorem2}, together assumption (A1), we have that $\big| {\cal Z}(\ebo)\big| \rightarrow \infty$, in probability, as $n \rightarrow \infty$. 
\hspace*{\fill}$\blacksquare$

%%%%%%%%%%%%%%%%%%%%%%%%%%%%%%%%%%%%%%%
%%%%%%%%%%%%%%%%%%%%%%
\subsection{Lemmas}

\hh In order to prove Propositions \ref{proposition1}, \ref{corolary4}, \ref{proposition2}, \ref{proposition3} and Theorems \ref{theorem2} and \ref{theorem_H1}, we need the following lemmas. \\
The first four lemmas establish equivalent results obtained by \cite{Guo:Zou:Wang:Chen:13} for a linear model without change-point. \\
%%%%%%%%%%%%%%%%%%%%%%%%%%%%%%%%%%%%%%%%%%%%%%%%%%%%%%%%%%%%%%%%%%%
%%%%%%%%%%%%%%%%%%%%%%%%%%%%   Lemma 1  %%%%%%%%%%%%%%%%%%%%%%%%%%%
%%%%%%%%%%%%%%%%%%%%%%%%%%%%%%%%%%%%%%%%%%%%%%%%%%%%%%%%%%%%%%%%%%%
Let us consider the following $p \times p$ matrix $\textbf{L}_n \equiv \es^0_n -\ev_n^0$. For $1 \leq u,v \leq p$ let us consider $L_{n;u,v}$ the $(u,v)$ element of the matrix $\textbf{L}_n$. Let also consider the largest absolute element of $\textbf{L}_n$ : $\max _{1 \leq u,v \leq p} | S^0_{n,(u,v)}-V^0_{n,(u,v)}|$, with $S^0_{n,(u,v)}$ and $V^0_{n,(u,v)}$ denote the $(u,v)$ components of $\es^0_n$ and $\ev^0_n$, respectively.
\begin{lemma} 
\label{lemma1}
Under null hypothesis $H_0$, suppose that assumptions (A3) and (A4) are satisfied. For any $\epsilon >0$, there exists a positive constant $C_q$ that depending only on $q \geq 4$, such that 
\begin{equation*}
\eP\big [ \max _{1 \leq u,v \leq p} | S^0_{n,(u,v)}- V^0_{n,(u,v)}| \geq \epsilon  \big ] \leq C_q \frac{ \,p^2}{n^{q/2} \epsilon^q}.
\end{equation*}
\end{lemma}
\textbf{Proof.} The matrix $\textbf{L}_n$ can be written
\begin{equation*}
\label{eq180}
%L_n&=& \bigg [  \frac{1}{n \theta ^2 } \sum  _ {i=1} ^ {n\theta}  \ez_i^0 \ez_i^{0t} + \frac{1}{n (1-\theta )^2} \sum _{j =n\theta+1} ^ n \ez_j^0 \ez_j^{0t}\bigg ]
% - \bigg [  \frac{1}{n \theta ^2 } \sum  _ {i=1} ^ {n\theta} \ev_{(i)}^0 +\frac{1}{n (1-\theta )^2} \sum _{j =n\theta+1} ^ n \ev_{(j)}^0 \bigg ]
%\nonumber \\ 
%&=& \bigg [  \frac{1}{n \theta ^2 } \sum  _ {i=1} ^ {n\theta}  \XX_i (\varepsilon_i \varepsilon_i^t) \XX_i^t + \frac{1}{n (1-\theta )^2} \sum _{j =n\theta+1} ^ n  \XX_j (\varepsilon_j \varepsilon_j^t) \XX_j^t \bigg ] 
%\nonumber \\ &&
%- \bigg  [  \frac{1}{n \theta ^2 } \sum  _ {i=1} ^ {n\theta} \Var (\XX_i) \Var(\varepsilon_i) +\frac{1}{n (1-\theta )^2} \sum _{j =n\theta+1} ^ n \Var (\XX_j) \Var(\varepsilon_j) \bigg ]
%\nonumber \\
\textbf{L}_n= \frac{1}{n \theta ^2 } \sum  _ {i=1} ^ {k}  \XX_i  \XX_i^t (\varepsilon_i^2-\sigma^2) + \frac{1}{n (1-\theta )^2} \sum _{j =k+1} ^ n  \XX_j \XX_j^t (\varepsilon_j^2 -\sigma^2).
\end{equation*}
The $(u,v)$-th element of the matrix $\textbf{L}_n$, for $1 \leq u,v \leq p$, is 
\begin{equation}
\label{EqLn}
L_{n,(u,v)}=  \frac{1}{n \theta ^2 } \sum  _ {i=1} ^ {k}  X_{i,u}  X_{i,v} (\varepsilon_i^2-\sigma^2) + \frac{1}{n (1-\theta )^2} \sum _{j =k+1} ^ n  X_{j,u}  X_{j,v} (\varepsilon_j^2 -\sigma^2)
\equiv L^{(1)}_{n,(u,v)} +L^{(2)}_{n,(u,v)}. 
\end{equation}
Since $\theta=k/n \rightarrow \theta^0 \in (0,1)$ as $n \rightarrow \infty$, we can apply Lemma 1 of \cite{Guo:Zou:Wang:Chen:13} for $ L^{(1)}_{n,(u,v)}$ and $ L^{(2)}_{n,(u,v)}$. Then,  for fixed $\theta \in (0,1)$, for all $\epsilon >0$ there exists two positive constants $C^{(1)}_q$ and $C^{(2)}_q$ such that:
\begin{equation}
\label{Lol}
\eP \big[| L^{(1)}_{n,(u,v)} | \geq \epsilon \big] \leq C^{(1)}_q \frac{p^2}{n^{q/2} \epsilon^q}.
\end{equation}
and 
\begin{equation}
\label{L2n2}
\eP \big[| L^{(2)}_{n,(u,v)} | \geq \epsilon \big] \leq C^{(2)}_q \frac{p^2}{n^{q/2} \epsilon^q}.
\end{equation}
Then the Lemma follows from relations (\ref{EqLn}), (\ref{Lol}) and (\ref{L2n2}), considering $C_q=\max\big(C^{(1)}_q,C^{(2)}_q\big)$.
\hspace*{\fill}$\blacksquare$  \\

By the next lemma we prove first that  all eigenvalues of $\es_n^0$ converge to those of $\ev_n^0$ uniformly with the rate $O_p(p\,\max _{1 \leq u,v \leq p} | S^0_{n,(u,v)}-V^0_{n,(u,v)}|)$ and then that all eigenvalues of $\es_n^0$ are bounded and strictly positive for $n$ enough large. 
%%%%%%%%%%%%%%%%%%%%%%%%%%%%%%%%%%%%%%%%%%%%%%%%%%%%%%%%%%%%%%%%%%%
%%%%%%%%%%%%%%%%%%%%%%%%%%%%   Lemma 2  %%%%%%%%%%%%%%%%%%%%%%%%%%%
%%%%%%%%%%%%%%%%%%%%%%%%%%%%%%%%%%%%%%%%%%%%%%%%%%%%%%%%%%%%%%%%%%%
\begin{lemma}
\label{corolary2}
Under assumptions (A1), (A3), (A4) and (A6), there exists two constants $C_0, C_1 >0$ such that the inequality $\gamma_p(\es_n^0) > C_0$ and $\gamma_1(\es_n^0) < C_1$ hold with probability tending to one as $n \rightarrow \infty$.
\end{lemma}
\textbf{Proof.}
Arguments based on the properties of the trace and of the eigenvalues of symmetric square matrices, together with  Lemma 2 of \cite{Guo:Zou:Wang:Chen:13} imply that 
\begin{equation}
\label{lemme2}\max  _{ 1 \leq r \leq p} |\gamma_r(\ev_n^0)-\gamma_r(\es_n^0)| \leq p\,\max _{1 \leq u,v \leq p} | S^0_{n,(u,v)}-V^0_{n,(u,v)}|.
\end{equation}
Lemma follows  by combining this last relation with  Lemma \ref{lemma1} and assumption (A1).
\hspace*{\fill}$\blacksquare$  \\ 

%%%%%%%%%%%%%%%%%%%%%%%%%%%%%%%%%%%%%%%%%%%%%%%%%%%%%%%%%%%%%%%%%%
%%%%%%%%%%%%%%%%%%%%%%%%%%% Lemme 3%%%%%%%%%%%%%%%%%%%%%%%%%%%%%%%%ùù
%%%%%%%%%%%%%%%%%%%%%%%%%%%%%%%%%%%%%%%%%%%%%%%%%%%%%%%%%%%%%%%%%%%%%%%%%%%%%ù
Let us consider the following random variable $ T_n^0 \equiv \max _{\substack{1\leqslant i \leqslant k\\  k+1 \leqslant j \leqslant n}} \left \{ \displaystyle \theta^{-1}  \|\ez_i ^0\|, (1-\theta)^{-1} \|\ez_j ^0\| \right \}$.
\begin{lemma}
\label{lemma3}
Under the null hypothesis $H_0$, suppose that assumptions (A3) and (A4) are satisfied. Then, for $q \geq 4$, we have 
$T_n^0=o_{\eP}(p^{1/2}n^{1/q} ). $
%where $T_n^0=\max   _{ 1 \leq i \leq n} \|\ez_i ^0||$.
\end{lemma}
\textbf{Proof.} 
Since $\theta=k/n \rightarrow \theta^0 \in (0,1)$ as $n \rightarrow \infty$, we have with probability 1, for enough large $n$: $T_n^0 \leq C \max _{ 1 \leq i \leq n} \|\ez_i ^0\| $, with $C>0$. By Lemma 3 of \cite{Guo:Zou:Wang:Chen:13}, we have that: $ \max _{ 1 \leq i \leq n} \|\ez_i ^0\| =o_{\eP}( p^{1/2}n^{1/q}) $ and the lemma follows.
\hspace*{\fill}$\blacksquare$  \\
%%%%%%%%%%%%%%%%%%%%%%%%%%%%%%%%%%%%%%%%%%%%%%%%%%%%%%%%%%%%%%%%%%%
%%%%%%%%%%%%%%%%%%%%%%%%%%%%   Lemma 4  %%%%%%%%%%%%%%%%%%%%%%%%%%%
%%%%%%%%%%%%%%%%%%%%%%%%%%%%%%%%%%%%%%%%%%%%%%%%%%%%%%%%%%%%%%%%%%%

By the following Lemma we give an asymptotic approximation for the $L_2$-norm of the vector $\epsi_n(\ebo)$, given by (\ref{psib}), under hypothesis $H_0$.
\begin{lemma} 
\label{lemma4}
Under the  null hypothesis $H_0$, if assumption (A1) holds, we have 
$\|\epsi_n^0\|= O_{\eP}( p^{1/2}n^{-1/2}).$
\end{lemma}
\textbf{Proof.} Let  $\XX^{(1)}$ is the $p \times k$ design matrix whose $k$ columns are $ \XX_i$, for $i=1,\ldots,k$ and $\XX^{(2)}$ is the $p \times (n-k)$ design matrix whose $(n-k)$ columns are $ \XX_j$, for $j=k+1,\ldots,n$. Since $(\varepsilon_i)$ are independent, we have that
\begin{eqnarray*}
\label{eq28}
\eE \big [\epsi_n^{0t}\epsi_n^0 \big ]
&=&\frac{1}{(n\theta)^2} \eE \big[  \sum  _{i,i'=1} ^{k}  \ez_i^{0t} \ez_{i'}^{0} \big ]
+ \frac{1}{(n(1-\theta))^2}\eE \big [  \sum  _{j,j'=k+1} ^{n} \ez_j^{0t}  \ez_{j'}^{0} \big ]
%\nonumber \\
%&=&\frac{\sigma^2}{n}  tr \Big( \frac{1}{n\theta^2} (\XX^{(1)})^t  \XX^{(1)}\Big ) +\frac{\sigma^2}{n}  tr \Big( \frac{1}{n(1-\theta)^2} (\XX^{(2)})^t  \XX^{(2)}  \Big)
%\nonumber \\
%&=&\frac{\sigma^2}{n}  tr \Big( \frac{1}{n\theta^2}  \XX^{(1)} (\XX^{(1)})^t  \Big )+\frac{\sigma^2}{n}  tr \Big( \frac{1}{n(1-\theta)^2}   \XX^{(2)} (\XX^{(2)})^t \Big)
\nonumber \\
&=&\frac{\sigma^2}{n}  \tr \Big( \frac{1}{n\theta^2}  \XX^{(1)} (\XX^{(1)})^t  +\frac{1}{n(1-\theta)^2}   \XX^{(2)} (\XX^{(2)})^t \Big).
\nonumber \\
\end{eqnarray*}
Using assumption (A1), we obtain that
\begin{eqnarray*}
\label{eq282}
\eE \big [\epsi_n^{0t}\epsi_n^0 \big ]
&\leq& \sigma^2 \frac{p}{n}  \bigg (  \gamma_1 \Big(\frac{1}{n \theta}   \XX^{(1)} (\XX^{(1)})^t + \frac{1}{n (1-\theta)} \XX^{(2)} (\XX^{(2)})^t \Big) \bigg)
=O(pn^{-1}).
\end{eqnarray*}
Then, 
$
\|\epsi_n^0\|=O_{\eP}(p^{1/2}n^{-1/2} )
$.
\hspace*{\fill}$\blacksquare$ \\

The following Lemma gives a first approximation for the EL statistic, under hypothesis $H_0$.
\begin{lemma}
\label{debut_Prop2}
Under the same assumptions as in Proposition \ref{proposition2}, if hypothesis $H_0$ is true, we have:
\[
\EL_{nk}(\ebo)= 2n \el ^t\epsi_n^0 -n\el^t \es_n^0 \el +\frac{2}{3} \bigg (  \frac{1}{\theta^3}  \sum  _ {i=1} ^ {k}  (\ez_i^{0t} \el)^3- \frac{1}{(1-\theta)^3} \sum _{j =k+1} ^ n  (\ez_j^{0t}\el)^3\bigg)  +o_{\eP}(1),
\]
with $\theta=k/n.$
\end{lemma}
\textbf{Proof.} 
 The limited development of the statistic $\EL_{nk}(\ebo)$ specified by  relation (\ref{eq6}), in the neighbourhood of $\el=\textbf{0}_p$, up to order 3 can be written
\begin{eqnarray}
\label{DL}
\EL_{nk}(\ebo)&=& \el^t\bigg ( \frac{2}{\theta}  \sum_ {i=1} ^ {k}  \ez_i^0 - \frac{2}{1-\theta } \sum _{j =k+1} ^ n \ez_j^0\bigg) - 
 \el^t\bigg ( \frac{1}{\theta^2}  \sum  _ {i=1} ^ {k}  \ez_i^0 \ez_i^{0t} + \frac{1}{(1-\theta )^2} \sum _{j =k+1} ^ n \ez_j^0\ez_j^{0t} \bigg) \el 
\nonumber  \\ &&
+ \frac{2}{3} \bigg(  \frac{1}{\theta^3}  \sum_ {i=1} ^ {k}  (\ez_i^{0t} \el)^3- \frac{1}{(1-\theta)^3} \sum _{j =k+1} ^ n  (\ez_j^{0t}\el)^3\bigg) 
%\nonumber  \\ &&
+ \frac{1}{4!}   \sum_ {u,v,r,s=1}^{p}  
\frac{\partial^4 \EL_{nk}(\ebo) (\widetilde{\el}_{uvrs}) }{\partial \lambda_u  \partial \lambda_v \partial \lambda_r \partial \lambda_s} (\lambda_u)(\lambda_v) (\lambda_r) (\lambda_s)  
\nonumber  \\ 
& \equiv & {\cal E}_1 -{\cal E}_2+{\cal E}_3+{\cal E}_4,
\end{eqnarray}
where for all $1 \leq u,v,r,s \leq p$, $\lambda_u$ is the u-th component of $\el$ and $\widetilde{\el}_{uvrs}=a_{uvrs} \el$, with $a_{uvrs} \in [0,1]$. \\ 
  We first study ${\cal E}_4$, which can be written
\begin{eqnarray}
\label{last}
{\cal E}_4&=& -\frac{1}{4!} \bigg (  \frac{12}{\theta^4}  \sum  _ {i=1} ^ {k}  \frac{(\ez_i^{0t}  \el)^4}{(1+\theta^{-1} \widetilde{\el} \ez_i^0  )^4}+\frac{12}{(1-\theta)^4} \sum _{j =k+1} ^ n  \frac{(\ez_j^{0t} \el)^4}{(1-(1-\theta)^{-1} \widetilde{\el} \ez_i^0  )^4}\bigg).
 \end{eqnarray}  
By  Proposition 1 of \cite{Ciuperca:Salloum:15}, we have that,  for all $\epsilon >0$, there exists two positive absolute constants $M_1$ and $M_2$ such that
\begin{equation*}
\eP \bigg [ \frac{1}{M_1 \theta^4 } \sum  _ {i=1} ^ {k}  (\ez_i^{0t} \el)^4 \leq   \sum  _ {i=1} ^ {k}  \frac{(\ez_i^{0t}  \el)^4}{(1+\theta^{-1} \widetilde{\el} \ez_i^0  )^4} \leq  \frac{1}{M_2 \theta^4 } \sum  _ {i=1} ^ {k}  (\ez_i^{0t} \el)^4 \bigg ] \geq 1- \epsilon.
\end{equation*}
Then, for the first term of the right-hand side of (\ref{last}), applying Cauchy-Schwartz's inequality, we obtain that 
\begin{equation*}
\sum  _ {i=1} ^ {k}  \frac{(\ez_i^{0t} \el)^4}{(1+ \theta^{-1} \widetilde{\el} \ez_i^0  )^4}   \leq    \frac{1}{M_2 \theta^4 } \sum  _ {i=1} ^ {k}  (\ez_i^{0t}  \el)^4
\leq \frac{1}{M_2 \theta^4 } \sum  _ {i=1} ^ {k}  \|\XX_i \|^4  \varepsilon_i^4 \| \el\|^4.
\end{equation*}
Using assumptions (A1) (A2), together with the fact that $\|\el\|=  O_{\eP}( p^{1/2}n^{-1/2})$ given by Proposition \ref{proposition1} and $p=o(n^{1/2})$, we obtain that the first term of (\ref{last}) is $o_{\eP}(1)$. In the same way we can demonstrate that each term of (\ref{last}) is $o_{\eP}(1)$, which implies that ${\cal E}_4=o_{\eP}(1)$.  \\  

Using notations given by (\ref{S0}) and (\ref{psi}), we obtain that ${\cal E}_1= 2 n \el^t \epsi_n^0 $ and ${\cal E}_2=  n \el^t \es_n^0 \el$. Then, relation (\ref{DL}) becomes
$
\EL_{nk}(\ebo)= 2n \el ^t\epsi_n^0 -n\el^t \es_n^0 \el +{\cal E}_3+o_{\eP}(1)
$
and the lemma follows.
\hspace*{\fill}$\blacksquare$ \\

%%%%%%%%%%%%%%%%%%%%%%%%%%%%%%%%%%%%%%%%%%%%%%%%%%%%%%%%%%%%%%%%%%%
%%%%%%%%%%%%%%%%%%%%%%%%%%%%   Lemma 5  %%%%%%%%%%%%%%%%%%%%%%%%%%%
%%%%%%%%%%%%%%%%%%%%%%%%%%%%%%%%%%%%%%%%%%%%%%%%%%%%%%%%%%%%%%%%%%%
The following Lemma gives, under hypothesis $H_0$,  an approximation for the Lagrange multiplier $\el$ and the asymptotic behaviour of $\sup_{ 1 \leq i \leq n }   |\el^t\ez_i^0|=o_{\eP}(1)$.
\begin{lemma}
\label{lemma5}
Under the null hypothesis $H_0$, suppose that assumptions (A1), (A3)-(A6) are satisfied. Then for $q \geq 4$ and fixed $\theta \in (0,1)$, we have $\sup_{ 1 \leq i \leq n }   |\el^t\ez_i^0|=o_{\eP}(1)$  and
\begin{equation}
\label{ettql}
\el=(\es_n^0)^{-1}(\epsi_n^0 + o_{\eP}(n^{(1-q)/q} p^{3/2})) .
\end{equation}
\end{lemma}
\textbf{Proof.} By Proposition \ref{proposition1}, we have that $\|\el\|= O_{\eP}(n^{-1/2} p^{1/2})$. Note that, by Cauchy-Schwartz's inequality and by Lemma \ref{lemma3}, we have
\[\theta^{-1} \sup  _{ 1 \leq i \leq k} |\el^t\ez_i^0| \leq \|\el \|  T_n^0=O_{\eP}(n^{-1/2} p^{1/2})  o_{\eP}(n^{1/q} p^{1/2})= o_{\eP}(n^{(-q+2)/2q} p )\]
and 
\[(1-\theta)^{-1}\sup  _{ k+1 \leq j \leq n} |\el^t\ez_j^0| \leq \|\el \|  T_n^0=O_{\eP}(n^{-1/2} p^{1/2}) \ o_{\eP}(n^{1/q} p^{1/2})= o_{\eP}( n^{(-q+2)/2q} p ),\]
with $\theta=k/n$.
These two relations together assumption (A5) involve $\sup_{ 1 \leq i \leq n }   |\el^t\ez_i^0|=o_{\eP}(1)$. \\ 
We prove now  relation (\ref{ettql}). 
The limited development of (\ref{eq120}), in the neighbourhood of $\el=\textbf{0}_p$, up to order 3 can be written
\begin{equation}
\begin{array}{lll}
\label{eq200}
 \textbf{0}_p&=&\displaystyle{ \Big (\frac{1}{\theta}  \sum  _ {i=1} ^ {k}  \ez_i^0-\frac{1}{1-\theta}  \sum  _ {j=k+1} ^ {n}  \ez_j^0 \Big )  
  - \Big (\frac{1}{\theta^2}  \sum  _ {i=1} ^ {k}  \ez_i^0 \ez_i^{0t} +  \frac{1}{(1-\theta)^2}  \sum  _ {j=k+1} ^ {n} \ez_j^0 \ez_j^{0t} \Big )  \el }
  \\ &&
   +\displaystyle{ \Big(\frac{1}{\theta^3}  \sum  _ {i=1} ^ {k}  \ez_i^0  ( \ez_i^{0t}  \el )^2-  \frac{1}{(1-\theta)^3} \sum  _ {j=k+1} ^ {n}  \ez_j^0 ( \ez_j^{0t}  \el )^2 \Big ) }  \\ &&
+\displaystyle{ \Big(\frac{1}{(\theta + {\widetilde \el} ^{t}  \ez_i^0)^4}  \sum  _ {i=1} ^ {k}  \ez_i^0  ( \ez_i^{0t}  \el )^3-  \frac{1}{(1-\theta - \widetilde { \el} ^{t}  \ez_j^0)^4}  \sum  _ {j=k+1} ^ {n}  \ez_j^0 ( \ez_j^{0t}  \el )^3 \Big ),}
\end{array}
\end{equation}
where $ \widetilde {\el} = u \el$, with $ u \in (0,1)$. \\ 
Using Proposition 1 of \cite{Ciuperca:Salloum:15}, similarly as for the term ${\cal E}_4$ of Proposition \ref{proposition2}, we can demonstrate easily that the last term of the right hand side of relation (\ref{eq200}) is $o_{\eP}(1)$.
We recall that
\begin{equation}
\label{R1n}
\er_n^0  \equiv  \frac{1}{n\theta^3}  \sum  _ {i=1} ^ {k}  \ez_i^{0}  ( \el ^t \ez_i^0)^2 - \frac{1}{n(1-\theta)^3} \sum  _ {j =k+1} ^ n   \ez_j^{0}  ( \el ^t \ez_j^0)^2.
\end{equation}
Then, using notations given by (\ref{S0}), (\ref{psi}) and (\ref{R1n}) we obtain that  relation (\ref{eq200}) becomes $\epsi_n^0 -\es_n^0 \el +\textbf{R}_n ^0=\textbf{0}_p$. Thus 
\begin{equation}
\label{eq21}
\el= (\es_n^0)^{-1}(\textbf{R}_n^0+ \epsi_n^0)(1+o_{\eP}(1)).
\end{equation}
We recall that $ T_n^0 \equiv \max _{\substack{1\leqslant i \leqslant k\\  k+1 \leqslant j \leqslant n}} \left \{ \displaystyle \theta^{-1}  \|\ez_i ^0\|, (1-\theta)^{-1} \|\ez_j ^0\| \right \}$. 
Then, we have for  $\er_n^0$:
\[
\big \| \er_n^0\big \|
 \leq    T_n^0    \el ^t \bigg( \frac{1}{n\theta^2}   \sum  _ {i=1} ^ {k} \ez_i^0 \ez_i^{0t} +\frac{1}{n(1-\theta)^2}   \sum  _ {j=k+1} ^ {n}   \ez_j^0 \ez_j^{0t} \bigg ) \el
 \leq   T_n^0 \el ^t  \es_n^0 \el.
\]
Using Lemma 4 of \cite{Liu:Zou:Wang:13} and Lemma \ref{corolary2}  we obtain that 
$
\big \| \er_n^0\big \| \leq  \ T_n^0    \|\el\|^2 \gamma_1(\es_n^0)=O_{\eP}( T_n^0    \|\el\|^2)$.
 On the other hand, by Proposition \ref{proposition1} we have that  $\|\el\|= O_{\eP}(n^{-1/2} p^{1/2} )$ and  by Lemma \ref{lemma3} that  $T_n^0=o_{\eP}(n^{1/q} p^{1/2})$. Then for $\|\er_n^0 \|$, we obtain that
\begin{eqnarray}
\label{eq23}
\|\er_n^0 \| 
%&= & O_{\eP}( T_n^0   \,  ||\el||^2) 
%\nonumber \\ 
%= O_p(n^{-1} p)\, o_{\eP}(n^{1/q} p^{1/2} )
=o_{\eP}( n^{(1-q)/q}  p^{3/2} ).
\end{eqnarray}
Relation (\ref{ettql}) follows from (\ref{eq21}) and (\ref{eq23}).
\hspace*{\fill}$\blacksquare$  \\

%%%%%%%%%%%%%%%%%%%%%%%%%%%%%%%%%%%%%%%%%%%%%%%%%%%%%%%%%%%%%%%%%%%
%%%%%%%%%%%%%%%%%%%%%%%%%%%%   Lemma 6  %%%%%%%%%%%%%%%%%%%%%%%%%%%
%%%%%%%%%%%%%%%%%%%%%%%%%%%%%%%%%%%%%%%%%%%%%%%%%%%%%%%%%%%%%%%%%%%
The following result gives an asymptotic  approximation for $\big[(\es_n^0)^{-1} -(\ev_n^0)^{-1}\big] \epsi_n^0$ and $\big[(\es_n^0)^{-1} -(\ev_n^0)^{-1}\big] \er_n^0$. 
\begin{lemma}
\label{lemma6}
Under the null hypothesis ($H_0$), if assumptions (A3), (A4) and (A6) hold, we have \\ 
(i) $\big((\es_n^0)^{-1} -(\ev_n^0)^{-1}\big) \epsi_n^0  = \big((\ev_n^0)^{-1} \epsi_n^0 \big) o_{\eP}(1).$ \\
(ii) $\big((\es_n^0)^{-1} -(\ev_n^0)^{-1}\big) \er_n^0  =\big((\ev_n^0)^{-1} \er_n^0\big) o_{\eP}(1).$
\end{lemma}
\textbf{Proof.} (i) 
By Lemma \ref{lemma1} we have that, under assumptions (A3) and (A4), for all $\epsilon >0$  there exist $C_q>0$, for $q \geq 4$, such that 
\begin{equation*}
\label{eqas}
\eP\big [p\ \max _{1 \leq u,v \leq p} | S^0_{n,(u,v)}-V^0_{n,(u,v)}| \geq \epsilon  \big ] \leq C_q \frac{ p^{2+q}}{n^{q/2} \epsilon^q}.
\end{equation*}
Furthermore,  under assumption (A6) we have
\begin{equation}
\label{remark2}
\max _{1 \leq u,v \leq p} | S^0_{n,(u,v)}-V^0_{n,(u,v)}| =o_{\eP}(1).
\end{equation} 
We recall that, for a matrix \textbf{A}, $\| \textbf{A} \|_1$ is the subordinate norm to the vector norm $\| . \|_1$.   Using Lemma 1$(iii)$ of \cite{Ciuperca:Salloum:15}, Lemma \ref{corolary2},  relations  (\ref{lemme2}),  (\ref{remark2}) and the identity $\ev_n^0((\es_n^0)^{-1} -(\ev_n^0)^{-1})\epsi_n ^0 
=(\ev_n^0 (\es_n^0)^{-1} -\textbf{I}_p )\epsi_n^0 $, we have  
\begin{eqnarray*}
\label{eq30} 
\big\| \ev_n^0 \big((\es_n^0)^{-1} -(\ev_n^0)^{-1}\big) \epsi_n^0 \big\| 
&=& \big\|\big (\ev_n^0 -\es_n^0\big )  (\es_n^0)^{-1} \epsi_n ^0  \big\|
\leq \big \|\ev_n^0 -\es_n^0  \big \|_1  \big \|(\es_n^0)^{-1} \big \|_1   \big \| \epsi_n^0  \big\|
\nonumber \\ 
& \leq & \max _{1\leq r \leq p} |\gamma_r(\ev_n^0 -\es_n^0) |\cdot |\gamma_1(\es_n^0)^{-1}| \cdot \big \|\epsi_n^0  \big \|
  \leq  C p\max _{1 \leq u,v \leq p} | S^0_{n,(u,v)}-V^0_{n,(u,v)}|  \cdot \|\epsi_n^0  \|
\nonumber \\ 
&=& \big \| \epsi_n ^0 \big\| o_{\eP}(1),
\end{eqnarray*}
which implies that
%$\ev_n^0 \big[(\es_n^0)^{-1} -(\ev_n^0)^{-1} \big] \epsi_n^0 =\epsi_n^0  $. Then 
$\big((\es_n^0)^{-1} -(\ev_n^0)^{-1}\big) \epsi_n^0=(\ev_n^0)^{-1} \epsi_n^0  o_{\eP}(1)$. 
 \\ 
(ii) The proof of (ii) is similar to (i).  
\hspace*{\fill}$\blacksquare$  \\ 

%%%%%%%%%%%%%%%%%%%%%%%%%%%%%%%%%%%%%%%%%%%%%%%%%%%%%%%%%%%%%%%%%%%
%%%%%%%%%%%%%%%%%%%%%%%%%%%%   Lemma 7  %%%%%%%%%%%%%%%%%%%%%%%%%%%
%%%%%%%%%%%%%%%%%%%%%%%%%%%%%%%%%%%%%%%%%%%%%%%%%%%%%%%%%%%%%%%%%%%
The following lemma is needed for proving  Proposition \ref{proposition3}.
We recall that $\ek_n^0= \textbf{I}_p-(\ev_n^0)^{-1/2} \es_n^0 (\ev_n^0)^{-1/2} $.
\begin{lemma}
\label{lemma7}
Under null hypothesis $H_0$, if assumption (A7) holds, then
\[\tr(\ek_n^0)^2=O_{\eP}(p^2 n^{-1}). \]
\end{lemma}
\textbf{Proof.} We  show that $\eE[\tr(\ek_n^0)^2]=O(p^2 n^{-1}).$  
For this, we write the  matrix $\ek_n^0$ as
\begin{equation*}
\label{eq110}
\ek_n^0 = \textbf{I}_p  -  \frac{1}{n \theta ^2 } \sum  _ {i=1} ^ {k}  \ew_i^0 \ew_i^{0t}  \frac{1}{n (1-\theta )^2} \sum _{j =k+1} ^ n \ew_j^0 \ew_j^{0t}. 
\end{equation*}
Then
\begin{equation}
\label{eq112}
\begin{array}{lll}
\eE[\tr (\ek_n^0)^2 ] &= &n^{-2}\sum _{r,s=1} ^ p \eE \big [  \big( \theta ^{-2 } \sum  _ {i=1} ^ {k}  w_{i,r}^0  w_{i,s}^0 + \frac{1}{(1-\theta )^2} \sum _{j =k+1} ^ n w_{j,r}^0 w_{j,s}^0 \big)^2  \big]  \\
 & & \\
& & -2 n^{-1}  \sum _{r=1} ^ p \big( \theta ^{-2 } \sum  _ {i=1} ^ {k}  \eE [w_{i,r}^0 w_{i,r}^0 ]+ (1-\theta  )^{-2} \sum _{j =k+1} ^ n \eE[w_{j,r}^0 w_{j,r}^0]\big) + p. 
\end{array}
\end{equation}
For the first term of the right-hand side of (\ref{eq112}), using the independence of $\ew_i^0$ for all $i=1,\cdots,n$, we have that 
\[
\frac{1}{n^2} \sum _{r,s=1} ^ p \eE \big[  \big(\frac{1}{ \theta ^2 } \sum  _ {i=1} ^ {k}  w_{i,r}^0 w_{is}^0+ \frac{1}{ (1-\theta )^2} \sum _{j =k+1} ^ n w_{j,r}^0 w_{j,s}^0 \big)^2\big]
\] 
 \[ \quad =\frac{1}{n^2}  \sum _{r,s=1} ^ p  \big (  \frac{1}{ \theta ^4 } \sum  _ {i=1} ^ {k}  \eE[ w_{i,r}^0 w_{i,r}^0 w_{i,s}^0 w_{i,s}^0 ]+ 
\frac{1}{ (1-\theta )^4} \sum _{j =k+1} ^ n \eE[ w_{j,r}^0 w_{j,r}^0 w_{j,s}^0 w_{j,s}^0] \big)
\]
\[\qquad \qquad +\frac{1}{n^2}  \sum _{r,s=1} ^ p \big (  \frac{1}{\theta ^2 } \sum  _ {i=1} ^ {k}  \eE [w_{i,r}^0  w_{i,s}^0] + \frac{1}{ (1-\theta )^2} \sum _{j =k+1} ^ n \eE [w_{j,r}^0 w_{j,s}^0] \big)^2
\]
\begin{equation}
\label{eq113}
\quad = \frac{1}{n}  \sum _{r,s=1} ^ p \alpha^{rrss}  +  \sum _{r,s=1} ^ p (\alpha^{rs} )^2 
=O(p^2 n^{-1}).
\end{equation}
The last equation is due to assumption (A7)  and to the fact that $\alpha^{rs} =0$, for $r\neq s$. 
\\ 
For the second term of the right-hand side of (\ref{eq112}), by the fact that $\alpha^{rr} =1$, we have that 
\begin{equation}
\label{eq114}
 \frac{1}{n} \sum _{r=1} ^ p \bigg( \frac{1}{ \theta ^2 } \sum  _ {i=1} ^ {k}  \eE [w_{i,r}^0 w_{i,r}^0 ]+ \frac{1}{ (1-\theta )^2} \sum _{j =k+1} ^ n \eE[w_{j,r}^0 w_{j,r}^0]\bigg)
=   \sum _{r=1} ^ p \alpha^{rr} 
= p.
\end{equation}
Then, using relations (\ref{eq112}), (\ref{eq113}) and (\ref{eq114}), we obtain that
\begin{equation*}
\label{eq115}
\eE[\tr (\ek_n^2 )]=O(p^2 n^{-1}).
\end{equation*}
By Markov's inequality, Lemma yields.
\hspace*{\fill}$\blacksquare$ \\ \\
 
%\section*{References}
 
%\bibliographystyle{apalike}
%\textbf{\textsl{References}}


\begin{thebibliography}{plain}



%%%%%%%%%%%%%%%%%%
\bibitem[Chen et al.(2009)]{Chen:Peng:Qin:09}
Chen, S. X., Peng, L., Qin, Y. L.,
 \newblock {Effect of data dimension on empirical likelihood}. \newblock {\textit{Biometrika}},
 \textbf{96}, 1-12, 2009.
%%%%%%%%%%%%%%%%%
\bibitem[Chow and Teicher(1997)]{Chow:Teicher:97}
Chow, Y.S., Teicher, H.,
\newblock {Probability theory: independence, interchangeability, martingales}. \newblock {\textit{Springer}},  New York, 1997.
  %%%%%%%%%%%%%%%%%
  \bibitem[Ciuperca(2014)]{Ciuperca:13a}
Ciuperca, G.,  
\newblock Model selection by LASSO methods in a change-point model, 
\newblock{\textit{Statistical Papers}}, \textbf{55}, 349-374, 2014.
\bibitem[Ciuperca(2013)]{Ciuperca:13b}
Ciuperca, G.,  
\newblock {Quantile regression in high-dimension with breaking}, 
\newblock{\textit{Journal of Statistical Theory and Applications}}, \textbf{12}, 288-305, 2013.
  %%%%%%%%%%%%%%%%% 
 \bibitem[Ciuperca and Salloum(2015)] {Ciuperca:Salloum:15}
Ciuperca, G., Salloum, Z., 
 \newblock {Empirical likelihood test in a posteriori change-point
nonlinear model}.
 \newblock To appear,  \textit{Metrika},  DOI 10.1007/s00184-015-0534-z.
 %%%%%%%%%%%%%%%%%  

 %%%%%%%%%%%%%%%%%
\bibitem[Dicker et al.(2013)]{Dicker:Huang:Lin:13}
Dicker, L., Huang, B., Lin, X., 
\newblock Variable selection and estimation with the seamless-$L_0$ penalty.
\newblock {\it Statistica Sinica}, \textbf{23},  929-962, 2013.
%%%%%%%%%%%%%%%%%%
\bibitem[Fan and Peng (2004)]{Fan:Peng:04}
Fan, J., Peng, H., 
\newblock Nonconcave penalized likelihood with a diverging number of parameters.
\newblock {\it The Annals of Statistics}, \textbf{32}(3),  928-961, 2004.
\bibitem[Guo et al(2013)]{Guo:Zou:Wang:Chen:13}
Guo, H., Zou, C., Wang, Z., Chen, B.,
 \newblock {Empirical likelihood for high-dimensional linear regression models}. \newblock {\textit{Metrika}},
 \textbf{77}(7), 921-945, 2013.
% %%%%%%%%%%%%%%%%%
 %%%%%%%%%%%%%%%%%  
%%%%%%%%%%%%%%%%%%
\bibitem[Hjort et al.(2009)]{Hjort:Mckeague:Vankeilegom:09}
Hjort, N.L., Mckeague, I.W., Van Keilegom, I.,
 \newblock {Extending the scope of empirical likelihood}. \newblock {\textit{The Annals of Statistics}},
 \textbf{37}, 1079-1111, 2009.
% %%%%%%%%%%%%%%%%%
\bibitem[Huang et al.(2012)]{Huang:Breheny:Ma:12}
Huang,  J., Breheny, P., Ma, S.,
 \newblock {A selective review of group selection in high-dimensional models}. \newblock {\textit{Statistical Science}},
 \textbf{27}(4), 481-499, 2012.
 \bibitem[Lee et al.(2015)]{Lee:Seo:Shin:12}
Lee S., Seo M.H., Shin Y., 
\newblock The LASSO for high-dimensional regression with a possible change-point.
\newblock \textit{Journal of the Royal Statistical Society: Series B}, DOI: 10.1111/rssb.12108.
\bibitem[Liu et al.(2008)]{Liu:Zou:Zhang:08}
Liu, Y.,  Zou,  C.,  Zhang, R.,  
\newblock {Empirical likelihood ratio test for a change-point in linear regression model}. \newblock {\textit{Communications in Statistics-Theory and Methods}},  \textbf{37}, 2551-2563,  2008. 
%%%%%%%%%%%%%%%%%%
\bibitem[Liu et al.(2013)]{Liu:Zou:Wang:13}
Liu, Y., Zou, C., Wang, Z.,
 \newblock {Calibration of the empirical likelihood for high-dimensional data}. \newblock {\textit{Annals of the Institute of Statistical Mathematics}}, \textbf{65}, 529-550, 2013.
 %%%%%%%%%%%%%%%%%  
 \bibitem[Lung-Yut-Fong et al.(2013)]{Lung-Yut-Fong:Levy:Cappe:12}
Lung-Yut-Fong, A., L\'evy-Leduc, C., Capp\'e, O.,
 \newblock {Distributed detection/localization of change-points in high-dimensional network traffic data},
 \newblock {\textit{Statistics and Computing}}, \textbf{22}, 485-496, 2012.

%%%%%%%%%%%%%%%%%
 
 
%%%%%%%%%%%%%%%%%%
%\bibitem[Owen(1990)]{Owen:90}
%Owen, A. B.,
%\newblock {Empirical likelihood ratio confidence regions}. \newblock {\textit{The Annals of Statistics}},  \textbf{18}, 90-120,  1990.
%%%%%%%%%%%%%%%%%
\bibitem[Tibshirani(1996)]{Tibshirani:96}
Tibshirani, R.,
\newblock {Regression shrinkage and selection via the LASSO.}
\newblock{\it Journal of the Royal Statistical Society}, Ser. B, {\bf 58}, 267-288,1996.
%%%%%%%%%%%%%%%%%
\bibitem[Wu and Liu(2009)]{Wu:Liu:09}
Wu,  Y.,  Liu, Y.,
\newblock {Variable selection in quantile regression.}
\newblock {\textit{Statistica Sinica}}, \textbf{19}, 801-817, 2009.
%%%%%%%%%%%%%%%%%%
\bibitem [Zi et al.(2012)]{Zi:Zou:Liu:12}
Zi, X., Zou, C., Liu, Y.,
 \newblock {Two-sample empirical likelihood method for difference between coefficients in linear regression model}.
 \newblock {\textit{Statistical Papers}}, \textbf{53}(1), 83-93, 2012.
 \bibitem[Zou and Yuan(2008)]{Zou:Yuan:08}
Zou, H., Yuan, M.,  
\newblock {Composite quantile regression and the oracle model selection theory.}
\newblock{\textit{The Annals of Statistics}}, \textbf{36}, 1108-1126, 2008.

\end{thebibliography}
\end{document}